\def\draft{n}
\documentclass[12pt]{amsart}
\usepackage[pagewise]{lineno}
\usepackage[headings]{fullpage}
\usepackage{amssymb,epic,eepic,epsfig, amsmath}
\usepackage{pb-diagram,tikz, mathtools, enumitem}
\usetikzlibrary{matrix, arrows}

\theoremstyle{plain}

\newtheorem{theorem}{Theorem}
\newtheorem{construction}{Construction}
\newtheorem{proposition}{Proposition}[section]
\newtheorem{lemma}[proposition]{Lemma}
\newtheorem{corollary}[proposition]{Corollary}

\newtheorem*{lawrence*}{Lawrence's Theorem}
\newtheorem*{brion*}{Brion's Theorem}

\theoremstyle{definition}
\newtheorem{definition}[proposition]{Definition}

\theoremstyle{remark}
\newtheorem{example}[proposition]{Example}

\def\printname#1{
        \if\draft y
                \smash{\makebox[0pt]{\hspace{-0.5in}
                        \raisebox{8pt}{\tt\tiny #1}}}
        \fi
}

\newlength{\standardunitlength}
\catcode`\@=11
\long\def\@makecaption#1#2{%
     \vskip 10pt

\setbox\@tempboxa\hbox{
       \small\sf{\bfcaptionfont #1. }\ignorespaces #2}%
     \ifdim \wd\@tempboxa >\captionwidth {%
         \rightskip=\@captionmargin\leftskip=\@captionmargin
         \unhbox\@tempboxa\par}%
       \else
         \hbox to\hsize{\hfil\box\@tempboxa\hfil}%
     \fi}
\font\bfcaptionfont=cmssbx10 scaled \magstephalf
\newdimen\@captionmargin\@captionmargin=2\parindent
\newdimen\captionwidth\captionwidth=\hsize
\catcode`\@=12

\newcommand{\Span}{\operatorname{Span}}

\def\lbl#1{\label{#1}\printname{#1}}

\newcommand{\ba}{\begin{align*}}
\newcommand{\ea}{\end{align*}}

\def\BZ{\mathbb Z}
\def\BP{\mathbb P}

\def\BR{\mathbb R}
\def\BC{\mathbb C}

\def\E{\mathcal E}

\def\M{\mathcal M}

\def\calR{\mathcal R}

\def\Cone{\mathrm{Cone}}

\def\Vol{\mathrm{Vol}}
\def\Td{\mathrm{Td}}
\def\td{\mathrm{td}}
\def\Star{\mathrm{Star}}

\def\calR{\mathcal{R}}

\def\calC{\mathcal{C}}
\def\calM{\mathcal{M}}
\def\calD{\mathcal{D}}

\def\lin{\mathrm{lin}}

\def\Supp{\mathrm{Supp}}

\def\dim{\mathrm{dim}}
\def\lin{\mathrm{lin}}

\def\lin{\mathrm{lin}}

\def\lin{\mathrm{lin}}
\def\aff{\mathrm{aff}}

\def\Vert{{\rm{Vert}}}

\begin{document}

\title{An algebraic construction of Sum-Integral Interpolators}
\author{Benjamin Fischer and
Jamie Pommersheim}
\maketitle

\begin{abstract}
    This paper presents an algebraic construction of Euler-Maclaurin formulas for polytopes.  The formulas obtained generalize and unite the previous lattice point formulas \cite{PT, Mo}, and Euler-Maclaurin formulas \cite{BV1, GP}.  While the approach of this paper originates in the theory of toric varieties, and recovers results from \cite{FP} proved using toric geometry, the present paper is self-contained and does not rely on results from toric geometry.  We aim in particular to exhibit in a combinatorial way ingredients, such as Todd classes and cycle-level intersections in Chow rings, that first entered the theory of polytopes from algebraic geometry. 
\end{abstract}
\maketitle

\section{Introduction}

{\bf 1.1 Overview.}\ \ Suppose $V$ is a vector space $V\cong \BR^n$ and $M\cong \BZ^n$ is a lattice in $V$.  If $P$ is a polytope in $V$ whose vertices lie in $M$, one may ask if there is a nice formula for $|P\cap M|$, the number of lattice points that are contained in $P$. Since one expects the volume of $P$ to give an approximation for the number of lattice points, one might ask to express the number of lattice points in $P$ as a linear combination of the volumes of the faces $F$ of $P$.  McMullen \cite{McM} asked in particular if it is possible to make the coefficients in this expression rational numbers that can be computed {\it locally} in the sense that the coefficient at each face $F$ depends only on the geometry of the polytope near $F$, specifically the tangent cone to $P$ at $F$.  McMullen showed nonconstructively that this is indeed possible, and decades later Morelli \cite{Mo} gave a construction at the cost of allowing coefficients that are not rational numbers, but rather live in a rational function field on a Grassmannian. The first construction with rational coefficents was achieved by Pommersheim and Thomas, [PT], who used the theory of toric varieties.  The construction of [PT] shows that any choice of a complement map on $V$, which is a certain choice of linear subspaces of $V$, yields a local formula of McMullen type. 

More general than the question of enumerating lattice points is the {\it Euler-Maclaurin question}: Given a function $f(x)$  on the space $V$, determine the sum of $f(x)$ over all lattice points $x\in P\cap M $. Berline and Vergne [BV] gave such an Euler-Maclaurin formula for summing polynomial functions $f(x)$ over rational polytopes $P$.  The formula of \cite{BV1} is defined inductively and relies on the choice of an inner product on $V$. In \cite{GP}, this construction was extended to general complement maps. The paper \cite{FP} contains a different construction of the \cite{BV1} coefficient $\mu$ in the case of an integral polytope, and contains a proof of the conjecture from \cite{GP} that the coefficients of the Berline-Vergne formula specialize to the earlier lattice point enumeration formula of Pommersheim-Thomas.  As in \cite{GP}, this construction extends to a general complement map. This algebraic construction of \cite{FP} and its proof are based in the theory of toric varieties. Guo, Paycha, and Zhang \cite{GPZ} reconstructed the Berline-Vergne formula for inner products using an Algebraic Birkhoff Factorization procedure, from the theory of Hopf algebras and non-commutative geometry \cite{CK}.

A main goal of this paper is to better understand the formulas arising from the construction of \cite{FP} and to obtain results of this construction in a simpler way that does not rely on the results or language from the theory of toric varieties. That is, we show how the choice of an inner product, or more generally a complement map, on a vector space $V$ can be used to give an algebraic construction of an Euler-Maclaurin formula for integral polytopes in $V$. We also give an explicit expression for the cone invariants appearing in these formulas. Along the way, we see how ingredients that entered the theory of polytopes from algebraic geometry, such as the Todd class and Chow rings, can also be seen naturally in a combinatorial way.  Many of these results are contained in a somewhat different and more general form in the first author's thesis \cite{Fi}.  

In the final section, we analyze the case in which the polytope is the standard simplex.  We show that a certain choice of complement map yields an equivariant generalization of the Fulton-Diaconis formula for the Todd class of projective space.

\bigskip
{\bf 1.2 Sum-Integral Interpolators.}\ \
Let $V\cong \BR^n$ be a real vector space, and $M\cong \BZ^n$ a full-dimensional lattice in $V$.  
Let $W$ be the dual space of all linear functionals on $V$. We will use $N$ to denote the lattice in $W$ dual to $M$.


There are two important valuations  (\cite{La,KP}), which we denote  $S$ and $I$ (for sum and integral)  on polyhedra $P$ in $V$.  The valuations $S$ and $I$ have values in the space $\M(W)$ of meromorphic functions on the complexification $W_{\BC}=W\otimes \BC$, and are continuations of the formulas, for $y\in W$, 
$$
S(P)(y) = \sum_{x\in P\cap M} e^{-\langle y, x \rangle},\qquad
I(P)(y) = \int_{x\in P} e^{-\langle y, x \rangle} \ dm_P,
$$
where volume $dm_P$ is normalized with respect to the sublattice of $M$ in the affine span of $P$, so that a fundamental domain has volume $1$.

It is well-understood that the problem of finding Euler-Maclaurin formulas for polytopes $P$ in $M$ reduces to being able to express $S(P)$ as a linear combination of $I(F)$, ranging over the faces $F<P$. (We use the notation $F<P$ to indicate that $F$ is a face of $P$.)  If we desire a local formula, then for each face $F$, we require the coefficient of $I(F)$ to depend only on the tangent cone $T(P,F)$, or equivalently, its dual $C(P,F)$.  That is, we would like to have a function $\mu$ on cones in $W$ such that for any integral polytope $P$ in $V$,
\begin{equation}\lbl{eq.si}
S(P) = \sum_{F<P} \mu(C(P,F)) I(F),
\end{equation}
where $\mu(C(P,F))$ lives in the ring $\calM_0(W)$ of meromorphic functions on $W_{\BC}$ that are regular at $0$. A function $\mu$ with this property is called an $SI${\it -interpolator}.  As shown in [BV], such a $\mu$ leads in a natural way to an Euler-Maclaurin formula for rational polytopes in $M$.   Specializing Equation (\ref{eq.si}) to $y=0$, one obtains a formula for the number of lattice points in $P$ in terms of the volumes of the faces $F$ of $P$.

Our main purpose here is to give a straightforward algebraic construction of  $\mu$.  In the end, our construction results in the same $\mu$ as in [GP]. In the case when $V$ comes with a chosen inner product, our construction results in the same $\mu$ as in [BV].  This algebraic construction was originally developed using the theory of toric varieties [FP], though our presentation here, both statements and proofs, do not require any results from that theory. 

\bigskip
{\bf 1.3 Choosing a complement map or an inner product.}\ \
To construct a $\mu$, we will require the choice of a complement map.  Introduced to the subject by Hugh Thomas \cite{Th}, a complement map $\Psi$ is a certain choice of linear subspaces.  
   In particular, if $L$ is a cone in $W$ of dimension $k$, then $\Psi(L)$ is a required to be a subspace of $V$, also of dimension $k$, that is complementary to the $(n-k)$-dimensional space $\lin(L)^\perp$, the space perpendicular to the linear span of $L$.   Complement maps need not be defined on all cones $L$.  As a result, our $\mu$ will be defined for cones which satisfy a genericity condition with respect to $\Psi$. See Section \ref{sec.polybasics}.3 for the precise definition of complement maps, $\Psi$-generic cones, and polytopes. 
   
   One of the primary ways complement maps arise is from  an inner product. Indeed, if the space $W$ comes endowed with an inner product, then there is a natural identification of $W$ with the dual space $V$.  One then gets a natural complement map on $W$ by taking $\Psi(L)$ to be the linear span of $L$. For readers who are interested primarily in the special case of our construction arising from inner products, we have included Section \ref{sec.ip}, which details the construction in this case as well as pointing out some simplifications in the proof. 
   
   Another way to get a complement map on $W$ is to choose a complete flag in $W$.  See Section \ref{sec.polybasics}, Example \ref{ex.flagcm}.
  Complement maps arising from inner products were used to give the first rational-valued local formula for the number of lattice points in a polytope in \cite{PT}. Complement maps arising from flags were used in \cite{PT} to understand the formulas of Morelli \cite{Mo}. Though complement maps most commonly arise in one of these two ways, in the final section of this paper, we relate the formulas of Diaconis and Fulton \cite{DF} to a complement map coming neither from an inner product nor from a flag.

\bigskip
{\bf 1.4 The main theorem.}\ \ Having briefly discussed sum-integral interpolators and complement maps, we are now in position to state the main theorem.

\begin{theorem}\lbl{thm.cm}
Let $W$ be a real vector space of dimension $n$ equipped with a complement map $\Psi$. Let $N$ be an $n$-dimensional lattice in $W$. Let $V=W^*$ be the dual space, and let $M\subset V$ be the lattice dual to $N$. Let $\calC^{\Psi}$ denote the set of all $\Psi$-generic pointed cones in $N$.   Then Construction \ref{constr.cm} (detailed in Section 1.5) results in a function 
$$
\mu : \calC^{\Psi} \longrightarrow \calM_0(W)
$$
satisfying
\begin{enumerate}
\item
 For any $n$-dimensional integral $\Psi$-generic polytope $P$ in $M$,
\begin{equation*}
S(P) = \sum_{F<P} \mu(C(P,F)) I(F).
\end{equation*}
\item
The function $\mu$ is additive under subdivisions. That is, if a $\Psi$-generic cone $L$ in $N$ is written as the union of $\Psi$-generic cones $L_i, i = 1\dots, r$, of the same dimension which intersect along faces of smaller dimension, then
$$
\mu(L) = \sum_{i=1}^r \mu(L_i).
$$

\end{enumerate}
\end{theorem}

The corresponding theorem in the case where the vector space $W$ comes with an inner product rather than a complement map appears as Theorem 2 in Section \ref{sec.ip}.

\bigskip\bigskip
{\bf 1.5 The construction.}\ \ We now give the construction of the interpolator $\mu$ in Theorem 1.  Before making the precise statement, we provide some motivation and introduce some of the main ingredients in the construction.  Details of this discussion are deferred to Section \ref{sec.interpbasic}.

To construct $\mu$, Brion's theorem \cite{Br} will allow us to work locally on the tangent cones of the polytope.  In this spirit, let us begin with an $n$-dimensional cone $L=\Cone(w_1, \dots, w_n)$ in $N$.  For now, we will assume that $L$ is a {\it basic} cone, meaning that it can be generated by a subset of lattice basis.  Let $K$ denote the dual cone, generated by the dual basis $\{v_1, \dots, v_n\} $  of $M$.  Our goal is to  express $S(K)$ as a linear combination of the $I(F)$, where $F$ runs over the faces of $K$. We wish the coefficients of this linear combination to be, like  $S(K)$ and $I(F)$, meromorphic functions on $W$, and we wish them to be regular at $0$, implying that they are represented in the power series ring
$$\Lambda = [M]]=\mathbb{C}[[v_1, \dots, v_n]].$$

 We now apply the ideas of Khovanskii-Pukhlikov \cite{KP} to our cone.  They achieved an Euler-Maclaurin formula for polytopes  by perturbing the polytope, moving its faces parallel to themselves independently, and then applying differential operators to the exponential integral over these perturbed cones.  With this in mind, we introduce the ring
 $$
 R(L)=\Lambda[[D_1,\dots, D_n]],
 $$
whose elements may be viewed as differential operators and applied to the integral over the the perturbed cones. This results in an {\it evaluation map}  $\E:R(L)\rightarrow\Lambda$.

On a basic cone, the sum $S(K)$ is simply a geometric series, which can be  obtained by applying the evaluation map $\E$ to the {\it Todd power series} $\td(D_1, \dots D_n)$, defined by
$$
\td(z_1, \dots, z_n) = \prod_{i=1}^n \frac{z_i}{1-e^{-z_i}}.
$$

The integral $I(F)$, where $F$ is a face of $K$, is also easy to compute; $I(F)$ can be obtained by applying $\E$ to a squarefree monomial in the $D_i$ that corresponds to $F$.

These simple observations have an interesting consequence: we can solve our problem of expressing $S(K)$ as a linear combination of the $I(F)$  if we can express $ \td(D_1, \dots D_n)$ as a  {\it squarefree} polynomial in the $D_i$, with coefficients in $\Lambda$.  Furthermore, we can work modulo the kernel of $\E$ in finding such an expression. 

However, there are two problems. First, such squarefree expressions are far from unique as the kernel of $\E$ is quite large.  In particular, for any $v\in V$, if we set
$$
l_v =  \sum_{i=1}^n \langle w_i, v \rangle  D_i, 
$$
then 
$
l_v - v \in \ker\E.
$
The second issue is that we would like to obtain a local formula, so that in our squarefree expression, the coefficient of each monomial depends only on the corresponding face of $L$, and not on the entire cone $L$.

It turns out that these non-uniqueness and locality issues can be solved simultaneously by replacing the kernel of $\E$ with a smaller ideal $J$ defined in Construction 1 below. In the ring $Z(L)=R(L)/J$, every element has a unique expression that is squarefree in the $D_i$, as stated formally in Proposition \ref{prop.unisfcm} in Section \ref{sec.cm}. Expressing  $ \td(D_1, \dots, D_n)$ in a squarefree manner, and applying $\E$ gives us the desired expression for $S(K)$ as a linear combination of $I(F)$.  

To extend this to non-basic cones, we use the fact that every cone $L$ in $N$ can be subdivided into a collection $\{ L_i\}$ of basic cones of the same dimension as $L$. We then define $\mu(L) = \sum \mu(L_i)$.  We will prove that this sum is independent of the chosen subdivision.

\begin{construction}\lbl{constr.cm} 
Let $W$ be a real vector space of dimension $n$ equipped with a complement map $\Psi$. Let $V=W^*$ be the dual vector space. Suppose that $N$ is an $n$-dimensional lattice in $W$.  Let $L=\Cone(w_1, \dots, w_k)$ be a $\Psi$-generic $k$-dimensional basic cone in $N$.  Given this data, we make the following definitions:
\begin{enumerate}[label=\alph*.]
\item Let $\Lambda= \mathbb{C}[[V]]$, the power series ring $\mathbb{C}[[v_1, \dots, v_n]]$ where $\{ v_1, \dots, v_n\} $ is any basis of $V$.
 \item Define the ring $R(L)$ as the power series ring over $\Lambda$ with indeterminates $D_1, \dots, D_k$ corresponding to the rays of $L$:
 $$
 R(L) = \Lambda[[D_1, \dots, D_k]]= \mathbb{C}[[v_1, \dots, v_n, D_1, \dots, D_k ]].
 $$
 \item For any face $G$ of $L$, denote by $D_G$ the corresponding squarefree monomial $\prod{D_i}$, with the product taken over all $i\in \{1, \dots, k \} $ such that $w_i$ is a ray of $G$. 
 \item For any $v\in V$, define a linear polynomial $l_v$ by
 $$
 l_v = \sum_{i=1}^k \langle w_i, v \rangle D_i \in R(L).
 $$
 \item Define an ideal $J(L)$ of $R(L)$ by
 $$
J(L) = \Biggl\langle D_G (l_{v} - v) \ \Biggl\rvert\ G<L, \ \ v\in \Psi(G)    \Biggr\rangle.
$$
\item Let 
$$
Z(L) = R(L)/ J(L).
$$
\item In $Z(L)$, let
$$
\Td(L) = \td(D_1, \dots, D_k) = \prod_{i=1}^{k} \frac{D_i}{1-e^{-D_i}}.
$$
\item   $\Td(L)$ has a unique expression in $Z(L)$ as a $\Lambda$-linear combination of squarefree monomials by Proposition \ref{prop.unisfcm}:
$$
\Td(L) = \sum_{G<L} \lambda_G D_G, \text{\ \ \ \ with\ \ } \lambda_G \in \Lambda.
$$
\item Define $\mu(L)$ as the coefficient of $D_1 \dots D_k$ in this expression, i.e., 
$$
\mu(L) = \lambda_L.
$$
\end{enumerate}

This defines $\mu(L)$ for $\Psi$-generic basic cones. If $L$ is any $\Psi$-generic cone of dimension $k$, then we define $\mu(L)$ by subdividing $L$ into $\Psi$-generic basic cones $L_i$ of dimension $k$, and setting
$$
\mu(L) = \sum \mu(L_i).
$$
Theorem \ref{thm.cm} includes the assertion that this is well defined. 
\end{construction}

The corresponding results for the special  case in which  $W$ comes with an inner product, rather than a general complement map, are spelled out in Construction 2 in Section \ref{sec.ip}.

Construction 1 is a direct generalization of the construction of \cite{PT}, which produces the constant terms $\mu_0(L):=\mu(L)(0)$, thereby giving a local formula for the number of lattice points in any polytope with vertices in the lattice $M$. 

Construction 1
 results in the same $\mu$ as in \cite{GP}, as proved in \cite{FP}. When $\Psi$ arises from an inner product, our construction results in the same $\mu$ as in \cite{BV1}.  Note, however that the formulas of \cite{BV1} and \cite{GP} apply to rational polytopes, whereas Constructions 1 and 2 require the polytope to be integral. The authors have so far been unable to extend the construction given here to the rational case.

\bigskip
{\bf 1.6\ Relation to toric varieties.}
It follows from the results in \cite{FP} that the coefficients arising out of Construction 1 express the equivariant Todd class of a toric variety.  Precisely, suppose $X_{\Sigma}$  is a  toric variety arising from a fan $\Sigma$ in $N$.   Given any complement map (or inner product) such that all cones of $\Sigma$ are $\Psi$-generic, we may compute  the values $\mu(L)$ for all cones of $\Sigma$.  It then follows that the equivariant Todd class has the following expression in terms of the orbit closures $V(L)$.
$$
\Td^T(X_{\Sigma}) = \sum_{L\in \Sigma} \mu(L) [V(L)].
$$

We also note the relation between the ring $R(L)$ and the equivariant cohomology ring of the toric variety.  The equivariant cohomology of a smooth toric variety $X_{\Sigma}$ has a presentation as a quotient of a power series ring $\Lambda[[D_1, \dots, D_s]]$ with variables corresponding to the rays of the fan (torus invariant divisors) (\cite{Br2}, see also \cite{Fu2}.)  The ideal of relations is a sum $I_1 + I_2$, where $I_1$ is the Stanley-Reisner ideal and $I_2$ consists of linear relations $l_v - v$, where $v$ is in the lattice dual to the fan.  Theorem 2 of \cite{FP} lifts the multiplication in this ring to the cycle level by replacing $I_2$ with a subideal $J^{\Psi}$ that depends on the complement map. That is, the ring $Z(\Sigma)$ obtained in this way is isomorphic to the cycle group, and its multiplication is a lifting of the multiplication on the equivariant cohomology ring.  Due to the presence of the Stanley-Riesner ideal, products in this ring may also be computed locally.  That is,  if $L$ is any cone of $\Sigma$, then coefficient of $[V(L)]$ in any product can be computed in the ring $Z(L)$ obtained from $Z(\Sigma)$ by modding out by all $D_i$ outside of $L$.  The presentation of the ring $Z(L)$ obtained in this way is exactly the one introduced in Construction 1.  

\bigskip
{\bf Acknowledgements.}   The authors would like to thank Federico Castillo, Sinai Robins, Quang-Nhat Le, and Evan Griggs for useful conversations. We also thank the anonymous referees for their useful comments.

\section{Basic Notions: Polytopes, Cones, and Complement Maps}\lbl{sec.polybasics}

In this section we introduce some well-known notions from the theory of polyhedra that we will need in this paper.  Section 2.1 is devoted to polytopes and cones. In Section 2.2, we discuss the properties of the valuations $S$ and $I$, including Brion's theorems.  Finally, Section 2.3 is devoted to complement maps, including definitions and basic examples.

\bigskip
{\bf 2.1\ Polytopes and cones.} Let $V\cong \BR^n$ be a real vector space, and $M\cong \BZ^n$ a full-dimensional lattice in $V$.  
Let $W$ be the vector space dual space to $V$.  Let $N$ be the dual lattice of $M$. This is the set of all linear functionals $w$ on $V$ that take integer values on all points of the lattice $M$:
$$
N= \{ w\in W \ |\ \langle w, v \rangle \in \BZ \ {\rm for\ all}\ v\in M \} 
$$
Here $\langle \cdot ,\cdot \rangle:W\times V\rightarrow \BR$ represents the natural pairing, evaluation of a functional on a vector.

A {\it cone} $K$ in $V$ is a subset of $V$ that is the convex hull of a finite number of rays from the origin.   For $v_1, \dots, v_k\in V$, the cone generated by these vectors is
$$
\Cone(v_1, \dots, v_k) = \biggl\{ \sum_{i=1}^k c_i v_i\  \biggl\rvert\  c_i \in \BR_{\geq 0} \biggr\}.
$$ 
All cones considered in this paper will be assumed to be {\it rational}, meaning that they can be generated by vectors $v_1, \dots, v_k$ that live in the lattice. A cone is called {\it pointed} if it does not contain a line.

If $K$ is a cone in $V$, then the {\it dual cone} of $K$ is the cone $L$ in $W$ defined by  
$$
L = \{ w\in W\ |\ \langle w, v \rangle \geq 0\ {\rm for\ all} \ v\in K \}.
$$

A {\it polyhedron} $P$ in $V$ is a subset of $V$ defined by a finite set of linear inequalities. The affine span of a polyhedron will be denoted $\aff(P)$. The vertex set of $P$ will be denoted $\Vert(P)$. All of the polyhedra we consider will be {\it integral}, meaning that $\Vert(P)$ is a subset of the lattice $M$. A {\it polytope} is a bounded polyhedron. Equivalently, a polytope is the convex hull of a finite set of points (its vertex set).

If $K$ is a cone in $V$, the linear span of $K$ will be denoted $\lin(K)$.  The perpendicular space $K^\perp$ is a subspace of $W$, defined as all linear functionals that vanish on $\lin(K)$.

If $P$ is a polyhedron in $V$ and $F$ is a face of $P$, then the {\it tangent cone} $T(P,F)$ is the cone of all directions in $V$ the one can move from $F$ and stay in $V$. Formally,
$$
T(P,F) = \{v \in V\ |\ x + \epsilon v \in P \rm{\ for\ small\ } \epsilon > 0 \}
$$ 
where $x$ is any point in the relative interior of $F$.

If $P$ is a polyhedron and $v$ is a vertex of $P$, then the {\it supporting cone} to $P$ at $v$ is the tangent cone $T(P,v)$ translated back to $v$. Formally, 
$
P_v= v+ T(P,v).
$

We define $C(P,F)$, the {\it normal cone to}\ $P$ {\it at}\ $F$, as the cone in $W$ dual to the tangent cone $T(P,F)$.  This is the cone in $W$ generated by the inner normals to the facets of $P$ that contain $F$.     As duals, the normal cone and the tangent cone contain the same information, but for convenience, we choose to work mainly with normal cones in this paper.

If $L$ is a cone in $W$ with dual cone $K$ in $V$, then there is a one-to-one order-reversing correspondence between the faces of $L$ and the faces of $K$.  In particular, for any face $G$ of $L$, the corresponding face of $K$ is  $K\cap G^\perp$.  Conversely, for a face $F$ of $K$, the corresponding face of $L$ is the normal cone $C(K,F)$.  This and other basic facts about dual cones can be found in [Ful], for example.

\bigskip
{\bf 2.2\ The valuations $S$ and $I$.}\ In this section, we lay out the important properties of the valuations $S$ and $I$, the exponential sum and integral, introduced in  (\cite{La,KP}), which we will need later.
For a rational polyhedron $P$ in $V$, these functions are given by the equations
\begin{equation}\lbl{eq.sandi}
    S(P)(y) = \sum_{x\in P\cap M} e^{-\langle y, x \rangle},\qquad
I(P)(y) = \int_{x\in P} e^{-\langle y, x \rangle} \ dm_P.
\end{equation}
for $y\in W_{\BC}$ whenever $|e^{-\langle y, x \rangle}|$ is summable (resp. 
integrable) over $P$.  The measure $dm_P$ denotes the relative Lebesgue measure on  $\aff(P)$  normalized with respect to the lattice $M\cap \aff(P)$. 

If $P$ is a polyhedron that does not contain a line, one sees that there will be values of $y$ for which the sum (resp. integral) converges absolutely, and as a result, the sum (resp. integral) has a unique meromorphic continuation to $W_{\BC}$. One extends $S$ and $I$ to all rational polyhedra by defining $S(P)=I(P)=0$ for any polyhedron $P$ that contains a line. Remarkably, it then turns out that $S$ and $I$ are {\it valuations}, meaning that they preserve relations among indicator functions of polyhedra.  The indicator function of $P$, which takes the value 1 on $P$ and 0 outside $P$, will be denoted by $[P]$.  

The valuation property of $S$ asserts that if we have a finite collection $\{P_i\}$ of polyhedra whose indicator functions satisfy
$\sum_i r_i [(P_i)]=0$, then the functions $S(P_i)$ satisfy the corresponding 
 relation 
$\sum_i r_i S(P_i)=0$.

The valuation property of $I$ is slightly different. Again suppose that we have a finite collection $\{P_i\}$ of polyhedra whose indicator functions satisfy
$\sum_i r_i [(P_i)]=0$. In this case,  the functions $I(P_i)$ satisfy 
$\sum_i r_i I(P_i)=0$,  with this latter sum restricted to the sum to those $P_i$ that do 
not lie in a proper affine subspace of $V$.  This is expressed by saying that $I$ is a {\it solid valuation}.

\begin{lawrence*}
Let $V$ be a finite-dimensional real vector space with $M$ a full-rank lattice in $V$, and let  $W$ denote the dual space. Then one can associate to each rational polyhedron  $P$ in $V$  two meromorphic functions $S(P)$,  the {\em exponential sum}, and $I(P)$, the {\em exponential integral} $I(P)$.  Both $S(P)$ and $I(P)$ live in the ring
 $\M(W)$ of meromorphic functions on the complexification $W_{\BC}=W\otimes \BC$.
The functions $S$ and $I$ satisfy the following properties.

\begin{itemize}
\item
If $P$ contains a straight line, then $S(P)=I(P)=0$.
\item
$S$  is a valuation.
\item
$I$ is a solid valuation. 
\item The values of $S(P)(y)$ and $I(P)(y)$ are given by Equation (\ref{eq.sandi}) whenever the sum or integral converges.

\end{itemize}

\end{lawrence*}

 Brion's Theorem \cite{Br,La} expresses $S(P)$, for any rational polyhedron $P$, in terms of $S$ evaluated at the vertex supporting cones of $P$. A similar relation holds for $I$.
\begin{brion*}  Let P be a rational polyhedron in V.  For each vertex $v$ of $P$, recall that $P_v$ is the supporting cone to $P$ at $v$. Then
\begin{align*}
    S(P) &= \sum_{v\in\Vert(P)}  S(P_v), {\rm \ \ \ and}\\
      I(P) &= \sum_{v\in\Vert(P)} I(P_v). 
\end{align*}
\end{brion*}

\bigskip
{\bf 2.3\ Complement Maps.}\ \ A complement map is a certain choice of complementary linear subspaces.   In particular, if $L$ is a cone in $W$ of dimension $k$, then we  assign to $L$ a subspace $\Psi(L)$ of $V$, also of dimension $k$.  We require that our choice for $\Psi(L)$ be generic in the sense that $\Psi(L)$ is complementary to the $(n-k)$-dimensional subspace $L^{\perp}$. 

\begin{definition} 
Let $W$ be a real vector space with dual space $V$.  A {\it complement map} on $W$ is a map $\Psi$ from a collection $\calD$ of cones of $W$ to subspaces of $V$ such that for any $L\in\calD$, 
\begin{equation}
    V= \Psi(L) \oplus L^{\perp},
\end{equation} 
and $\Psi$ preserves inclusion on its domain.  That is for any $U, U'\in\calD$, if $U\subset U'$, then $\Psi(U) \subset \Psi(U')$.

A  basic cone $L$ in $W$ is called {\it $\Psi$-generic} if for any face  $G$ of $L$, we have $G\in\calD$.   A non-basic cone $L$ is $\Psi$-generic if $L$ can be subdivided into $\Psi$-generic basic cones. A polytope $P$  in $V$ is $\Psi$-generic if for any face $F$ of $P$, the normal cone $C(P,F)$ is $\Psi$-generic.

\end{definition}

The above definition differs slightly from that of Thomas \cite{Th}, who introduced complement maps to the subject. 
\begin{example}\lbl{ex.ipcm}
If $W$ is a vector space, then the choice of an inner product on $W$, or equivalently the choice of an inner product on $V$, gives rise to a complement map.  For any cone $L$ in $W$, define $\Psi(L)\subset V$ to be the orthogonal complement (under the inner product) of the subspace $L^{\perp}\subset V$.  It is easy to check that this $\Psi$ is a complement map whose domain includes all cones in $W$. 
An equivalent way to say this takes advantage of the fact that the choice of an inner product on $W$ naturally identifies $W$ with its dual $V$. Under this identification, we simply set $\Psi(L)=\lin(L)$. 
\end{example}

\begin{example}\lbl{ex.flagcm}
If $W$ is a vector space, then the choice of a complete flag in $V$, or equivalently the choice of a complete flag in $W$, gives rise to a complement map. Indeed, let $\{0\}=V_0\subset V_1\subset \cdots \subset V_n=V$ be a complete flag in $V$, with $V_i$  being an $i$-dimensional subspace of $V$. Then
$$
\Psi(L) := V_{\dim(L)}
$$
defines a complement map whose domain consists of cones $L$ in $W$ such that $L^{\perp}$ is generic with respect to the flag, in the sense that $L^{\perp}$ meets $V_{\dim(L)}$ trivially.
\end{example}

\section{The interpolator property on basic cones}\lbl{sec.interpbasic}

In this section we introduce the evaluation map $\E$, which provides  motivation behind Construction \ref{constr.cm}, and use $\E$ to prove the interpolator property for basic cones.

Throughout, $V$ will be an $n$-dimensional real vector space with dual space $W$.  Inside $V$, let $M$ be a full-rank lattice, with $N\subset W$ being the lattice dual to $M$.

Let $K=\Cone(v_1,\dots, v_n)$ be an $n$-dimensional basic cone in $V$ with dual cone $L=\Cone(w_1,\dots, w_n)$.  Let
$$\Lambda = \mathbb{C}[[M]]=\mathbb{C}[[v_1, \dots, v_n]],$$ 
where $\{v_1, \dots, v_n\}$ is any basis of $M$. The ring $\Lambda$ is the completion of the polynomial algebra over $M$.

We compute 
$$
    I(K)(y) = \int_{x\in K} e^{-\langle y, x\rangle}\ dx =\frac1{\prod_{i=1}^n \langle y, v_i\rangle}.
$$ 
Notationally, since any $v\in V$ is a function on $W$, we can express this equation more briefly as
$$
    I(K) =\frac1{\prod_{i=1}^n  v_i},
$$ 
Here we view $I(K)$ as an element of the space  $\calM(W)$ of meromorphic functions on $W_{\BC}$. It may also be viewed as an element of  $\text{Frac}\ \Lambda$, the field of Laurent series in $v_1, \dots, v_n$.

 Following Khovanskii-Pukhlikov \cite{KP}, we introduce the
perturbed  cone $K_h$, defined for any $h=(h_1, \dots, h_n) \in \BR^n$ by
$$
K_h = \{ v\in V\ |\ \langle w_i, v \rangle \geq -h_i\}.
$$
An easy calculation shows that 
\begin{equation}\lbl{eq.ikhformula}
    I(K_h) = \frac1{\prod_{i=1}^n  v_i}\ e^{\sum h_i v_i}
\end{equation}

Let $R(L)$ denote the power series ring $\Lambda[[D_1,\dots, D_n]]$.  We view the elements of $R(L)$ as differential operators.

 We introduce the {\it evaluation map}, a $\Lambda$-linear map 
$$
\E : R(L) \rightarrow\text{Frac}\ \Lambda,$$ 
as follows. For any power series $p(D_1, \dots, D_n)$, we define 
$$
\E(p(D_1, \dots, D_n)) = p(\frac{\partial}{\partial h_1}, \dots, \frac{\partial}{\partial h_n})\diamond I(K_h)\ |_{h=0},
$$
where the $\diamond$ notation above indicates that the differential operator corresponding to the power series $p(D_1, \dots, D_n)$ is applied to the integral over the perturbed cone before being evaluated at $h=0$.
One calculates that 
\begin{equation}\lbl{eq.eeq}
    \E(p(D_1, \dots, D_n)) = \frac{1}{v_1 \cdots  v_n} p( v_1, \dots v_n)
\end{equation}
To show this, it is enough to check the case when $p$ is a monomial $p=D_1^{a_1}\cdots D_n^{a_n}$. In this case, using Equation (\ref{eq.ikhformula}),  we see that $ \E(p(D_1, \dots, D_n))$ becomes
$$
 (\frac{\partial}{\partial h_1})^{a_1}\cdots (\frac{\partial}{\partial h_n})^{a_n } \diamond \frac1{v_1\cdots v_n} e^{\sum h_i v_i}\ |_{h=0}= 
 \frac{1}{v_1 \cdots  v_n} v_1^{a_1}\cdots v_n^{a_n},
$$
verifying Equation (\ref{eq.eeq}).

If $G$ is a face of $L$, then we set 
$$
D_G := \prod D_i \in R(L),
$$
where the product is taken over all those $i$ such that $w_i$ is contained in $G$.  This establishes a one-to-one correspondence between the faces $G$ of $L$ and the squarefree monomials in $R(L)$.  


For any face $G$ of $L$, consider the corresponding face $K\cap G^{\perp}$ of $K$.  Equation (\ref{eq.eeq}) implies that
\begin{equation}\lbl{eq.sqfreee}
\E(D_G) = I(K\cap G^{\perp}),
\end{equation}
that is, the evaluation map $\E$ takes any squarefree monomial in $R(L)$ to the integral over the corresponding face of $K$.

\begin{example}\lbl{ex.2dcone} As an easy example, let $M=\BZ^2$, the integer lattice in  $V=\BR^2$.  
Let $K=\Cone(v_1, v_2)$ be the cone generated by the standard basis $\{v_1,v_2\}$ of $\BR^2$.   The perturbed cone, for $h=(h_1, h_2)\in \BR^2$, is given by
$$
K_h=\{(x_1, x_2) \in \BR^2 | x_1\geq -h_1, x_2\geq -h_2\}.
$$
We see that 
$$
I(K_h)(y) = \int_{x_1\geq -h_1, x_2\geq -h_2} e^{-\langle y,x\rangle} dx = \frac{e^{v_1h_1 + v_2 h_2}}{v_1 v_2},
$$
where, as above, $v_i$ is a shorthand for $\langle y, v_i\rangle$.
The evaluation map $\mathcal{E}$, the derivatives of this function at $h=0$, are given by 
$$
\E(D_1^{a_1} D_2^{a_2 })= (\frac{\partial}{\partial h_1})^{a_1} (\frac{\partial}{\partial h_2})^{a_2 } \diamond I(K_h)\ |_{h=0} =v_1^{a_1-1}v_2^{a_2-1}
$$
As an example of Equation (\ref{eq.sqfreee}), take $G=\Cone(w_1)$.  Then the dual face $K\cap G^{\perp}$ is $\Cone(v_2)$, so 
$$
I(K\cap G^{\perp})(y) = \int_{x_2\geq 0}  e^{-\langle y,(0, x_2)\rangle}\ dx= \frac1{v_2}, 
$$
which agrees with 
$\E(D_G) = \E(D_1) = \frac1{v_2}$.
\end{example}

Equation (\ref{eq.sqfreee}) ties the evaluation map $\E$ to the integrals $I(F)$ over the faces of $K$.
To tie the evaluation map $\E$ to the exponential sum $S(K)$, note that for our basic cone $K$, we have 
$$
S(K) = \sum_{v \in M \cap K} e^{-v} =\prod_{i=1}^n \frac{1}{1-e^{- v_i}}.
$$
Thus if we define  the {\it Todd power series} in $n$ variables by
$$
\td(z_1, \dots, z_n) = \prod_{i=1}^n \frac{z_i}{1-e^{-z_i}},
$$
then we see from Equation (\ref{eq.eeq}) that
\begin{equation}\lbl{eq.etds}
    \E( \td(D_1, \dots D_n)) = S(K).
\end{equation}

From Equations (\ref{eq.sqfreee}) and (\ref{eq.etds}), we see that whenever we have a squarefree expression for $ \td(D_1, \dots D_n)$ modulo $\ker \E$,  we get an  expression for $S(K)$ as a linear combination of the $I(F)$. This leads to the following proposition.

\begin{proposition}\lbl{prop.sibasic}
  Let $L$ be a basic $n$-dimensional cone in $W$ with dual cone $K$. Let $\E:R(L)\rightarrow\Lambda$ be the evaluation map defined above. Suppose that in $R(L)/\ker\E$ we have the squarefree expression
  \begin{equation}
    \td(D_1, \dots D_n)= \sum_{G<L} \lambda_G D_G, 
\end{equation}
 where the coefficients $\lambda_G$ are in $\Lambda$.  Then 
\begin{equation*}
S(K) = \sum_{G<L} \lambda_G\ I(K\cap G^{\perp}).
\end{equation*}
\end{proposition}

We note in particular that for any $v\in V$, if we set
$$
l_v =  \sum_{i=1}^n \langle w_i, v \rangle  D_i, 
$$
then by Equation (\ref{eq.eeq}),
$$
\E(l_v - v) =  \frac{1}{v_1 \cdots  v_n}\bigl(\sum_{i=1}^n \langle w_i, v \rangle  v_i - v\bigr) = 0.
$$
Thus, in Construction 1, when $L$ is an $n$-dimensional basic cone, the ideal $J(L)$ is indeed a subideal of $\ker\E$.  Thus the equation
$$
\Td(L) = \sum_{G<L} \mu(G) D_G,
$$
which holds in the ring $Z(L)=R(L)/J(L)$, implies that
\begin{equation}\lbl{eq.si2}
S(K) = \sum_{G<L} \mu(G) I(K\cap G^\perp).
\end{equation}
Using the one-to-one  correspondence between the faces of $L$ and the faces of $K$ described in Section 2.1, this is equivalent to the interpolator equation 
\begin{equation}\lbl{eq.siforcones}
S(K) = \sum_{F<K} \mu(C(K,F)) I(F),
\end{equation}
 for any basic $n$-dimensional cone $K$.

\section{Details and Proofs of the Complement Map Construction}\lbl{sec.cm}


In this section, we investigate the ring $Z(L)$ introduced in Construction \ref{constr.cm}.  Among other results, we prove the key fact that all elements of this ring can be expressed uniquely in a squarefree manner, as stated in   Proposition \ref{prop.unisfcm} below.

We adopt the notation of Construction 1. In particular, let $L$ be a $\Psi$-generic basic cone in $W$ generated by $\{w_1, \dots, w_k\}$.  As a notational convenience, in the rest of the paper we set
$
[k] = \{1,2,\dots,k\}
$, and for any $S\subset [k]$, we let 
$$
\Psi(S) := \Psi({\rm Cone}(\{w_s | s\in S\})).
$$
Thus the generators of $J(L)$ from Construction \ref{constr.cm} have the form $D_S(l_v-v)$, where $S\subset [k]$, and $v\in \Psi(S)$. 

\begin{lemma}\lbl{lem.standalone}  Let $L=\Cone(w_1,\dots,w_k)$ be a $\Psi$-generic cone.  Then
given $S\subset [k]$ and $i\in S$, there exists $v\in \Psi(S)$  such that $\langle w_i, v\rangle \neq 0$, but $\langle w_j,v \rangle  = 0 $ for all  $j\in S\setminus \{i\}$. 
\end{lemma}

\begin{proof}
Consider the map $f$ from $\Psi(S)$ to $\BR^{|S|}$ that sends $v$ to $f(v) = (\langle w_j,  v\rangle)_{j \in S}$.  We have
$$
\ker(f) = \Psi(S) \cap ({\rm span}\{w_s | s\in S\})^\perp =  \{0\}. 
$$
It follows that $f$ is injective, hence surjective, and the lemma follows. 
\end{proof}




\begin{proposition}\lbl{prop.unisfcm}  Let $W$ be an $n$-dimensional real vector space equipped with a complement map $\Psi$.  Let $L$ be a $\Psi$-generic basic cone in $W$ generated by the linearly independent set $\{w_1, \dots, w_k\}$ in $W$. Let  $Z(L)$ be the ring defined in Construction \ref{constr.cm}. 
  Then the squarefree monomials $\{D_G\ |\ G < L\}$ form a $\Lambda$-basis of $Z(L)$.
That is, every $\alpha\in Z(L)$ has a unique squarefree expression of the form
$$
\alpha = \sum_{G<L} \alpha_G D_G,
$$
with each $\alpha_G\in\Lambda$.
\end{proposition}

\begin{proof}   We first prove that any monomial $m=D_1^{a_1}\cdots D_k^{a_k}$ has a squarefree expression in $Z(L)$.   To do this, proceed by induction on the degree $\sum a_i$.  For degree 0 or 1, the monomial is already squarefree.  If our monomial is not squarefree, there is some $i$ such that $a_i\geq 2$.  By inductive assumption,  $m'=D_i^{-1} m$ has a squarefree expression.  Thus it is enough to show that for any squarefree monomial $D_S$, with $S\subset [k]$, the product $D_i D_S$ has a squarefree expression. If $i\notin S$, then we are done, so we can assume $i\in S$.

Let $v$ be as in Lemma \ref{lem.standalone}. Then in $Z(L)$, we have the relation $D_S(l_v-v)=0$, which means
\begin{equation}\lbl{eq.didscm}
\langle w_i,  v\rangle D_i D_S + \sum_{j\in [k]\setminus \{i\}} \langle w_j,  v\rangle D_j D_S- v D_S=0.
\end{equation}

Noting that all of the $ \langle w_j,  v\rangle$ vanish when $j\in S\setminus \{i\}$, we thus obtain the desired squarefree expression for  $D_i D_S$ in $Z(L)$.

To prove uniqueness, we must show that no squarefree polynomial belongs to $J(L)$ except $0$.  We fix $n$ and argue inductively on $k$, the dimension of $L$.  If $k=1$, then $\Psi(\{1\})$ is spanned by a single vector $u_1$, for which $\langle  w_1, u_1\rangle \neq 0$. Thus we see that the ring $R(L)$ is $\Lambda[[D_1]]$ and the ideal $J(L)$ is generated by $ D_1 (  \langle  w_1,u_1 \rangle D_1 - u_1 )$. This polynomial has no squarefree nonzero multiple except $0$.
        
Now suppose that $k>1$, and let $P$ be a squarefree polynomial in $J(L)$. In order to show that $P=0$, we will first show

\noindent {\bf Claim:} \ $P$ is divisible by $D_i$ for all $i=1,\dots, k$.

Without losing generality, we take $i=k$. Consider the cone $L'=\Cone(w_1, \dots w_{k-1})$. Since $L'$ is $\Psi$-generic, Construction \ref{constr.cm} yields the ring $R(L')$ containing the ideal $J(L')$, and by induction we may assume that $J(L')$ contains no  squarefree elements except $0$.   We will apply the homomorphism
$$
R(L)=\Lambda[[D_1,\dots, D_k]]\longrightarrow R(L')=\Lambda[[D_1,\dots, D_{k-1}]]
$$
that sends $D_k$ to $0$.   It is easy to see that under this map, the ideal $J(L)$ maps into $J(L')$. Indeed, for any $v\in V$, we see that
$$
l_{v}= \sum_{i=1}^{k} \langle w_i,v \rangle D_i\longmapsto l'_{v}= \sum_{i=1}^{k-1} \langle w_i,v \rangle D_i.
$$
Now take  $S\subset [k]$, and $v\in\Psi(S)$. Then if $k\in S$, the generator $D_S(l_v-v)$ maps to $0$, while if $k\notin S$, then the generator $D_S(l_v-v)$  maps to $D_S(l'_v-v)$, which is in $J(L')$. It follows that $J(L)$ maps to $J(L')$.

Since $P$ is squarefree in $J(L)$,  its image  $\bar{P}$ in $J(L')$ is also squarefree.  By inductive hypothesis, it follows that $\bar{P}=0$, which establishes the claim that $D_k$ divides $P$.

\vskip .25in
From the claim, it follows that $D_1\cdots D_k$ divides $P$.  Since $P$ is squarefree,   we must have $P=\alpha  D_1\cdots D_k$ for some $\alpha\in \Lambda$. Using Lemma \ref{lem.standalone}, we can find a basis $\{u_1, \dots, u_k\}$ of  $\Psi([k])$ such that $\langle  w_i,u_j \rangle =\delta_{ij}$.    Consider the $\Lambda$-linear homomorphism $e: R(L)\rightarrow \Lambda$ such that $e(D_i) = u_i$.   We show that $e(J)=0$. To see this, note first that $l_{u_i}= D_i$, so $e(l_{u_i}-u_i)=0$.  Since the $u_i$'s form a basis of $\Psi([k])$, this implies that $e(l_{v}-v)=0$ for any $v\in\Psi([k])$.  Since any generator of $J$ has a factor of the form $l_{v}-v$ with $v\in\Psi(G)$ for some face $G$, and $\Psi(G)\subset \Psi([k])$, we see that $e$ maps all generators of $J$ to 0.
As a result,  $P\in J$ implies $e(P) = \alpha u_1\cdots u_k = 0$.  It follows that  $\alpha=0$, and hence $P=0$, as desired.
 \end{proof}

Implicit in the proof above is an algorithm for finding the squarefree expression for any polynomial in $R(L)$, which amounts to repeatedly using Equation (\ref{eq.didscm}) with appropriately chosen $w$. 
One immediate and useful consequence of ({\ref{eq.didscm}}) is the following lemma, which asserts that when applying this algorithm, the support of a monomial can only grow.

\begin{lemma}\lbl{lem.suppgrow}
Given the notation of Proposition \ref{prop.unisfcm}, let $Q(D_1, \dots, D_k)$ be a monomial in $R(L)$.  Let $S=\Supp(Q)$ denote the support of $Q$, namely the set of variables that appear in $Q$ with positive exponent.  Then in the squarefree expression for $Q$ all nonzero terms have support containing $S$.
\end{lemma}

Another useful lemma says that each coefficient in a squarefree expression can be computed locally on the face $L'$ of $L$ corresponding to that monomial. This can also be seen from the algorithm itself.  Whenever a monomial appears with a variable that corresponds to a ray not in $L'$, we may ignore this monomial; by Lemma \ref{lem.suppgrow} the support will only grow, so the support of the squarefree expression of this monomial will not be contained in $L'$.  Here is the precise statement and proof.

\begin{lemma}\lbl{lem.facesf}
Suppose $W$, $L$, and $\Psi$ are as above.   Let $w_1, \dots, w_k$ be the generators of $L$ and consider any face $L'$ of $L$, which we may assume to be spanned by  \{$w_1, \dots, w_j\}$ with $j\leq k$,  Let $Q$ be any polynomial in $D_1, \dots, D_j$. 
If the squarefree expression for $Q$ in $Z(L)$ is 
$$
Q= \sum_{S\subset [k]} \alpha_S D_S,
$$
with $\alpha_S\in\Lambda$,
then the squarefree expression for $Q$ in $Z(L')$ is
$$
Q= \sum_{S\subset [j]} \alpha_S D_S,
$$
\end{lemma}

\begin{proof}
Take the equation in $R(L)$:
$$
Q= \sum_{S\subset [k]} \alpha_S D_S + \gamma,
$$
where $\gamma\in J(L)$.  Modulo $D_{j+1}, \dots, D_k$, generators of $J(L)$ become generators of $J(L')$ or $0$, and so one sees the equation becomes
$$
Q= \sum_{S\subset [j]} \alpha_S D_S + \gamma',
$$
with $\gamma'\in J(L')$. 
\end{proof}

We will also need to demonstrate the compatibility of the $\mu$ construction with inclusions. 
For this purpose, we examine how a complement map $W$ induces a complement map on any subspace $W'$. Note that if $W'$ is a subspace of $W$, then the dual to $W'$ is $V'=V/(W')^\perp$. For $v\in V$, we use $\overline{v}$ to denote the image of $v$ under the quotient map $V\rightarrow V'$, and use similar bar notation for the images of subspaces of $V$.

\begin{lemma}\lbl{lem.subpsi}
Let $W$ be a vector space equipped with a complement map $\Psi$.  Suppose that $W'$ is subspace of $W$. Define a complement map $\Psi'$ on $W'$ as follows.  For any cone $L$ in $W'$ in the domain of $\Psi$, set 
$$
\Psi'(L) = \overline{\Psi(L)},
$$
the image of $\Psi(L)$ under the natural projection $V\rightarrow V'$, where $V'$  is the vector space dual to $W'$.   Then $\Psi'$ is a complement map whose domain is all cones of $W'$ that are in the domain of $\Psi$.
\end{lemma}

\begin{proof}
Let $U=\lin(L)$. We must show that 
\begin{equation}\lbl{eq.compprime}
V'=\overline{\Psi(L)}\oplus \overline{U^\perp}
\end{equation}
Since $V=\Psi(L)\oplus U^\perp$, we need only check that the sum in (\ref{eq.compprime}) is direct, which follows from a dimension count: The kernel of the natural projection $V\rightarrow V'$ is $(W')^\perp$. Since  $\Psi(L)$ does not meet $U^{\perp}$, it is also the case that $\Psi(L)$ does not meet the kernel $(W')^\perp$. Thus, the dimension of $\Psi(L)$ remains unchanged under projection.  On the other hand, both $V$ and $U^\perp$ contain the entire kernel, so their dimensions fall by the same amount. Hence, the directness of $V=\Psi(L)\oplus U^\perp$ implies the directness in Equation (\ref{eq.compprime}).
\end{proof}

We now proceed to compatibility under inclusion.  As above, take $W'$ to be a subspace of $W$ and  $V'=V/(W')^\perp$ to be the dual to $W'$. Let $H$ be any subspace of $V$ that is complementary to $(W')^\perp$. Then the projection $V=H\oplus (W')^\perp\rightarrow H$ onto the first factor induces an isomorphism $V'\cong H$. Hence we obtain an inclusion $\iota: V'\rightarrow V$.  Note that for any $v\in H$, we have $\iota(\overline{v}) =v$.

Now suppose that $L$ is a basic cone in $W'$ in the domain of $\Psi$, and such that $\Psi(L)\subset H$. Considering $L$ as a cone in $W$, Construction \ref{constr.cm} yields $\Lambda=\BC[[V]], R(L), J(L), Z(L)$  and $\mu(L)$. On the other hand,  considering $L$ as a cone in $W'$, Construction 1 yields $\Lambda'=\BC[[V']], R'(L), J'(L), Z'(L)$ and $\mu'(L)$. The inclusion $\iota$ induces inclusions $\Lambda'\rightarrow \Lambda$ and  $R'(L)\rightarrow R(L)$, which we also denote $\iota$.

The following proposition asserts that $\mu$ is compatible with the inclusion map $\iota$. 
\begin{proposition}\lbl{prop.compatcm}
Let $W, W', \Psi,  \Psi'$ be as in Lemma \ref{lem.subpsi}. Let $H$ be any subspace of $V$ complementary to $(W')^\perp$.  Let $L=\Cone(w_1, \dots, w_k)$ be a basic cone  in $W'$ in the domain of $\Psi$ such that $\Psi(L)\subset H$. Let $\mu(L)$ (respectively $\mu'(L)$) be the result of applying Construction 1 to $L$ considered as a cone in $W$  (repectively $W'$). Let $\iota:\Lambda'\rightarrow \Lambda$ be the inclusion defined above, then we have
$$
\mu(L)=\iota(\mu'(L)).
$$
\end{proposition}

\begin{proof}
First, we prove that  under the inclusion $\iota: R'(L)\rightarrow R(L)$, we have $\iota(J'(L))\subset J(L)$. Note that the space $V'$ is dual to $W'$ with natural pairing $\langle w,\overline{v}\rangle = \langle w,v\rangle$ for any $v\in V$, and $w\in W'$. Thus for all $i$, we have $\langle w_i,\overline{v}\rangle = \langle w_i,v\rangle$, and it follows that $\iota(l_{\overline{v}})=l_{v}$.  Take any generator $D_S(l_{\overline{v}}-\overline{v})$, where $\overline{v}\in \Psi'(S) = \overline{\Psi(S)}$, with $v\in \Psi(S)$.  We see that $ v \in H$ and so $\iota(\overline{v})= v$. Hence the image of this generator under $\iota$ is exactly  $D_S(l_v-v)$, which is a generator of $J(L)$.   Thus $\iota$ maps $J'(L)$ into $J(L)$. 

By construction in $R'(L)$, we have
$$
\Td(L) = \sum_{G<L} \mu'(G) D_G +\gamma.
$$
with $\mu'(G) \in \Lambda'$ and $\gamma\in J'(L)$.  Applying $\iota$ yields
$$
\Td(L) = \sum_{G<L} \iota(\mu'(G)) D_G +\iota(\gamma).
$$
in the ring $R(L)$.
But then $\iota(\gamma)\in J(L)$, so the above equation gives a squarefree expression for $\Td(L)$ in $Z(L)$. Uniqueness of squarefree expressions in $Z(L)=R(L)/J(L)$ now forces $\mu(G)=\iota(\mu'(G))$ for all $G$.
\end{proof}


\section{Proof of Theorem \ref{thm.cm}}\lbl{sec.ipproof}

In this section, we complete the proof of Theorem \ref{thm.cm}.

Let $W$ be a real vector space of dimension $n$  equipped with a complement map $\Psi$, and let  $N$  be an $n$-dimensional lattice in $W$.  For all $\Psi$-generic basic cones $L$ in $W$, let $\mu(L)$ be as defined in Construction 1.

Recall that in Section \ref{sec.interpbasic}, we concluded that any $n$-dimensional basic cone satisfies the interpolator property of Equation (\ref{eq.si2}), or the equivalent formulation of Equation (\ref{eq.siforcones}).

It remains to show that 
 $\mu$ is well-defined, and additive under subdivisions, and satisfies the interpolator property for arbitrary integral polytopes.  To do this, we imagine the following set-up.  Let $L$ be an $n$-dimensional pointed cone in $W$, and suppose we have a subdivision of $L$. Precisely, let $\Sigma$ be a fan in $W$, a finite collection of cones closed under taking faces, and such that any two cones in $\Sigma$ intersect in a common face. We suppose that the union of the cones in $\Sigma$ is $L$, and we denote the $n$-dimensional cones in $\Sigma$ by $L_i$, $i=1,..., t$.   We use $K$ to denote the dual of $L$ and $K_i$ to denote the cone dual to $L_i$.

\begin{lemma}\lbl{prop.theprop}
Let $L$ be an $n$-dimensional pointed cone in $W$, and let  $\Sigma$ be a subdivision of $L$, as above. Fix a cone $\lambda\in \Sigma$, and let $\lambda_*$ denote the smallest face of $L$ containing $\lambda$.  Then we have the following congruence of indicator functions, modulo lines:
$$
\sum_{i:\lambda < L_i} [K_i \cap \lambda^\perp ] \equiv [K\cap \lambda_*^\perp]  {\rm\ \ \ \ modulo\ lines.}
$$
\end{lemma}

\begin{example}
Let $W$  be two-dimensional and let $L=\Cone(w_1, w_2)$ generated by a basis of $W$. (Nothing is lost here by picturing the standard basis of $\BR^2$.) Using the ray $\rho=\Cone(w_1+w_2)$, subdivide $L$ into two cones $L_1=\Cone(w_1, w_1+w_2)$ and  $L_2=\Cone(w_1+w_2, w_2)$.  Let $\Sigma$ be the resulting fan.  Denoting the dual basis by $\{v_1, v_2\}$, the dual cones are $K=\Cone(v_1, v_2)$, $K_1=\Cone(v_1-v_2, v_2)$ and  $K_2=\Cone(v_1, v_2-v_1)$.

Let us first take $\lambda$ to be the ray through $w_1$ and verify the conclusion of Lemma \ref{prop.theprop}. One finds $\lambda_*=\lambda$. The sum on the left side of the lemma has the single term $[K_1\cap \lambda^{\perp}]$, which equals $[\Cone(v_2)]$. On right hand side. $[K\cap \lambda_*^{\perp}]$ is also equal to  $\Cone(v_2)$.

Now take $\lambda$ to be the ray $\rho$ through $w_1+w_2$. Now the left hand side consists of the two terms $[\Cone(v_1-v_2)]$ and $[\Cone(v_2-v_1)]$.  The sum of these two rays is a line plus the origin. So modulo lines, the left hand side is simply $[\{0\}]$. Considering the right side, we see that  $\lambda_*=L$, so $[K\cap \lambda_*^\perp]=[\{0\}]$, and again the lemma checks out. 
\end{example}

Before proving Lemma \ref{prop.theprop}, we note the following corollary, which follows from the fact that $I$ is a solid valuation that vanishes on cones containing a line.  

\begin{corollary}\lbl{prop.icorolla}
With the hypotheses of Lemma \ref{prop.theprop},
$$
\sum_{i:\lambda < L_i} I(K_i \cap \lambda^\perp ) = 
\begin{cases} I(K\cap \lambda_*^\perp) & \text{if\ } \dim\ \lambda_* = \dim\ \lambda  \\ 
0 &\text{otherwise} 
\end{cases} 
$$
\end{corollary}

\begin{proof}[Proof of  Lemma \ref{prop.theprop}]
We consider the star $\Star\ \lambda$ (the smallest sub-fan of $\Sigma$ containing $\{ \sigma \in \Sigma\ | \lambda < \sigma \}$) and mod out by $U=\lin(\lambda)$, the linear span of $\lambda$. In this manner, in the quotient space $\overline{W}=W/U$, the fan $\Star\ \lambda$ defines a fan $\overline{\Star\  \lambda}$. 

{\bf Claim:\ } The union of the cones of $\overline{\Star\ \lambda}$ is  $\overline{L}$, the image of $L$ in $\overline{W}$.

{\it Proof of claim:} The fact that $\Star\ \lambda$ is contained in $L$ implies one inclusion. For the reverse inclusion, take  $w\in L$. Take a point $\rho$ in the relative interior of $\lambda$ and consider the path $f(t) = tw + (1-t)\rho$.  Note that $f(t)\in L$ for all $t\in[0,1]$. So there is some $L_i$ that contains $f(t)$ for all $t\in [0,\epsilon]$, for some $\epsilon>0$. That $L_i$ contains $f(0)=\rho$. Hence $L_i$ is in the star of $\lambda$. But then $f(\epsilon) = \epsilon w +(1-\epsilon) \rho\in L_i$ implies $\overline{w}\in \overline{L_i}$, as desired. This proves the claim.
\medskip

From the claim, we see that
$$
\sum_{i:\lambda < L_i} [\overline{L_i}] \equiv [\overline{L}] {\rm\ \ \ \  modulo\ cones\ of\ smaller\ dimension.}
$$
   
  We now dualize: The dual to $\bar{L}$ is $K\cap \lambda^\perp$ and the dual of  $\overline{L_i}$ is $K_i\cap \lambda^\perp$.  The dual to a smaller dimensional cone is a cone that contains a line. Since taking duals preserves linear relations among indicator functions, we see that 
$$
\sum_{i:\lambda < L_i} [K_i \cap \lambda^\perp ] \equiv [K\cap \lambda^\perp] {\rm\ \ \ \ modulo\ lines.}
$$
 The lemma now follows from the identity $K\cap \lambda^\perp= K\cap \lambda_*^\perp$. To see this, note that if $v\in K\cap\lambda^\perp$, then $G=L\cap v^\perp$ is a face of $L$ that contains $\lambda$, so $\lambda_*<G$.  Thus $v\in \lambda_*^\perp$.  The other inclusion is trivial.
\end{proof}

\bigskip
Having established the lemma and its corollary, we now prove that $\mu$ extends in a well-defined way to all pointed cones $L$. We proceed by induction on the dimension of $W$. Suppose that $W'$ is a proper subspace of $W$, and let $\Psi'$ be the complement map on $W'$ induced by $\Psi$, as in Lemma \ref{lem.subpsi}. The induction hypothesis tells us that Construction 1 yields  $\mu'$ that extends additively to all $\Psi'$-generic cones in $W'$. Since we have shown in Proposition \ref{prop.compatcm} that Construction 1 is compatible with the inclusion of $W'$ in $W$, it follows  that $\mu$ is well-defined and additive on  $\Psi$-generic cones of $W$ that have dimension smaller than the dimension of $W$. 

Let $L$ be a $\Psi$-generic cone in $W$ of dimenson $n=\dim\ W$, and let $\Sigma$ be any subdivision of $L$ into $\Psi$-generic basic cones $L_i$ of dimension $n$.  As usual, let $K$ and $K_i$ denote the respective dual cones. Then using Brion's Theorem, as well as Equation (\ref{eq.si2}), which has been established for $n$-dimensional basic cones, and Corollary \ref{prop.icorolla}, we have
\begin{align*}\lbl{eq.siadd}
    S(K) &= \sum_i S(K_i) \\
    &= \sum_i \sum_{\lambda:\ \lambda<L_i} \mu(\lambda) I(K_i\cap \lambda^\perp)\\
    &= \sum_{\lambda:\ \dim( \lambda_*)=\dim(\lambda)} \mu(\lambda) I(K\cap \lambda_*^\perp)\\
    &=\sum_{G<L} \biggl[ \sum_{\substack{\lambda:\ \lambda_*=G\\ \dim(\lambda_*)=\dim(\lambda)}}  \mu(\lambda) \biggr] I(K\cap G^\perp).
\end{align*}

Now consider the bracketed expression above.  For $G\neq L$, the bracketed expression equals exactly $\mu(G)$. To see this, note that the $\lambda$'s appearing in the bracketed sum form a collection of cones, all of the same dimension, that subdivide $G$.  Hence the induction hypothesis, the $\mu$ values of these cones add up to $\mu(G)$. Also by induction, all of these values $\mu(G)$ are independent of the subdivision.  For $G=L$, the bracketed expression equals $\sum \mu(L_i)$. This term is being multiplied by $I(K\cap L^\perp)=I(0) = 1$.  We thus have
\begin{equation}\lbl{eq.newsk}
     S(K) = \sum_{\substack{G<L\\ G\neq L}}   \mu(G)  I(K\cap G^\perp) + \sum_i \mu(L_i).
\end{equation}
Since $S(K)$ is independent of the subdivision, it follows that $\sum \mu(L_i)$ is also independent of the subdivision.  Hence $\mu(L)$ is well-defined.

It now follows that $\mu$ is an $SI$-interpolator for all $\Psi$-generic pointed cones: simply subdivide $L$ into $\Psi$-generic basic cones and use Equation (\ref{eq.newsk})  together with the definition $\mu(L)=\sum \mu(L_i)$. 

To prove that $\mu$ is additive under general subdivisions, again we may assume that this holds for smaller-dimensional cones. For a general subdivision of $L$  into $L_i$, the cones $K_i$ may no longer be basic.  However, we claim that the chain of equalities above, as well as Equation (\ref{eq.newsk}), still hold. Indeed, as we have just shown, the interpolator property is valid for each $K_i$, which justifies the second equality in the chain of equalities above, and all of the other equalities hold for the same reason as before. Now from Equation (\ref{eq.newsk}), and the fact that the interpolator property of Equation (\ref{eq.si2}) holds for $K$, if follows that  $\mu(L)=\sum \mu(L_i)$.

Next remark that if $K$ is a pointed cone translated by a lattice point $v\in M$, the $S$ and $I$ values are both multiplied by $e^{-v}$, while the $\mu$ values, which depend on normal cones to the faces, remain unchanged. Thus the interpolator equation holds in this case as well.

A standard argument using Brion's theorem  now shows that $\mu$ is an $SI$-interpolator on polytopes.  To wit, we have the following equalities, with justifications provided below.
\begin{align*}
    S(P) &= \sum_{v\in\Vert(P)} S(P_v)\\
    &= \sum_{v\in\Vert(P)} \sum_{F<P_v} \mu(C(P_v, F)) I(F)\\
    &= \sum_{v\in\Vert(P)} \sum_{E<P:v\in E} \mu(C(P,E)) I(E_v) \\
    &= \sum_{E<P} \mu(C(P,E)) \sum_{v\in\Vert(E)} I(E_v)\\
    &= \sum_{E<P} \mu(C(P,E)) I(E).
\end{align*}
The first and last equalities are Brion's theorem for $S$ and $I$, respectively. The second equality is the interpolator property for $P_v$. The fourth equality is interchange of summation.  Finally, the third equality follows from the fact that the faces $F$ of $P_v$ are in correspondence with the faces $E$ of $P$ that contain $v$, with $F=E_v$, the supporting cone to $E$ at $v$. Under this correspondence, we have an equality of normal cones  $C(P_v, F)= C(P,E)$.

Finally, we establish the meromorphicity of $\mu(L)$ and its regularity at 0. We use induction on $n$, the dimension of $V$.  If $n=0$, then $L$ is also $0$-dimensional, and $\mu(L)=1$. 

Now suppose $n>0$, and assume that $L$ is a basic cone in $N$.   If $\dim\ L <n$, then by the compatibility asserted in Proposition \ref{prop.compatcm}, $\mu(L)$ is the image of a value of $\mu$ from a smaller dimensional space, so we are done by induction.  Thus we can assume $\dim\ L =n$.  We then have, with the usual notation,
\begin{equation}
    S(K) = \sum_{G< L} \mu(G) I(K \cap G^{\perp}),
\end{equation}
 an equation taking place in $\frac1{v_1 \cdots v_n} \Lambda$.  That is, if we define
\begin{equation}
    \calR = \left\{ \rho \in \frac1{v_1 \cdots v_n} \Lambda \ \mid \ (v_1 \cdots v_n) \rho \text{\ defines a meromorphic function regular at 0.}\ \right\}, 
\end{equation}
then $S(K)\in\calR$, and $I(F)\in \calR$ for any $F<K$. By induction, we can assume that for any $G\neq L$,  $\mu(G)\in \calR$. Since $I(\{0\})=1$, it follows that $\mu(L)\in \calR$.  From this and the fact that $\mu(L)\in \Lambda$ is a power series, we conclude that $\mu(L)$ is in fact meromorphic and regular at the origin.  Since this holds for all basic cones, it holds for all pointed cones by additivity.

This completes the proof of Theorem \ref{thm.cm}.

\section{Details of the Inner Product Construction}\lbl{sec.ip}

An important special case of Construction 1 is the case in which the complement map $\Psi$ comes from an inner product.  In this section, we investigate this special case. We state the construction explicitly and mention a few ways in which this construction and its  proof are simpler than in the general case.

Notationally, we will use $[\cdot , \cdot]$ to refer to an inner product on $W$. (This distinguishes it from $\langle \cdot , \cdot \rangle$, which was used for a dual space acting on a vector space.)

 One immediate difference  is that  $\mu$ is now defined on all pointed cones. (Since $\Psi$ comes from an inner product, all  cones are $\Psi$-generic.) Here are the precise statements.

\begin{construction}\lbl{constr.ip} 
Let $W$ be a real vector space of dimension $n$ equipped with an inner product.  Suppose that $N$ is an $n$-dimensional lattice in $W$.  Let $L=\Cone(w_1, \dots, w_k)$ be a $k$-dimensional basic cone in $N$.  Given this data, we make the following definitions
\begin{enumerate}[label=\alph*.]
\item Let $\Lambda= \mathbb{C}[[W]]$. This is simply the power series ring $\mathbb{C}[[w_1, \dots, w_n]]$ where $\{ w_1, \dots, w_n\} $ is any basis of $W$.
 \item Define the ring $R(L)$ as the power series ring over $\Lambda$ with indeterminates $D_1, \dots, D_k$ corresponding to the rays of $L$:
 $$
 R(L) = \Lambda[[D_1, \dots, D_k]]= \mathbb{C}[[w_1, \dots, w_n, D_1, \dots, D_k ]].
 $$
 \item For any face $G<L$, denote by $D_G$ the corresponding squarefree monomial $\prod{D_i}$, with the product taken over all $i\in \{1, \dots, k \} $ such that $w_i$ is a ray of $G$. 
 \item For any $w\in W$, define a linear polynomial $l_w$ by
 $$
 l_w = \sum_{i=1}^k [w_i, w] D_i \in R(L).
 $$
 \item Define an ideal $J(L)$ of $R(L)$ by
 $$
J(L) = \Biggl\langle D_j (l_{w_j} - w_j) \ \Biggl\rvert\  j=1,\dots, k  \Biggr\rangle.
$$
\item Let 
$$
Z(L) = R(L)/ J(L).
$$
\item In $Z(L)$, let
$$
\Td(L) = \td(D_1, \dots, D_k) = \prod_{i=1}^{k} \frac{D_i}{1-e^{-D_i}}.
$$
\item   $\Td(L)$ has a unique expression in $Z(L)$ as $\Lambda$-linear combination of the squarefree monomials by Proposition \ref{prop.unisf}:
$$
\Td(L) = \sum_{G<L} \lambda_G D_G, \text{\ \ \ \ with\ \ } \lambda_G \in \Lambda.
$$
\item Define $\mu(L)$ as the coefficient of $D_1 \dots D_k$ in this expression, i.e., 
$$
\mu(L) = \lambda_L.
$$
\end{enumerate}

The above defines $\mu(L)$ for basic cones. If $L$ is any cone of dimension $k$, then we define $\mu(L)$ by subdividing $L$ into basic cones $L_i$ of dimension $k$, and setting
$$
\mu(L) = \sum \mu(L_i).
$$
Theorem \ref{thm.ip} (stated immediately below) includes the assertion that this is well defined. 
\end{construction}

\begin{theorem}\lbl{thm.ip}  Let $W$ be a real vector space of dimension $n$.  We assume that $W$ is equipped with an inner product, which allows identification of $W$ with the dual vector space $V$.   Suppose that $N$ is an $n$-dimensional lattice in $W$, and let $M$ denote the dual lattice. Then Construction \ref{constr.ip}  results in a well-defined function
$$
\mu : \{{\rm pointed\ cones\ in}\ N\} \longrightarrow \calM_0(W)
$$
satisfying
\begin{enumerate}
\item
 For any $n$-dimensional  polytope $P$ with vertices in $M$,
\begin{equation*}
S(P) = \sum_{F<P} \mu(C(P,F)) I(F).
\end{equation*}
\item
The function $\mu$ is additive under subdivisions. That is, if a cone $L$ in $N$ is written as the union of cones $L_i, i = 1\dots, r$ of the same dimension which intersect along faces of smaller dimension, then
$$
\mu(L) = \sum_{i=1}^r \mu(L_i).
$$

\end{enumerate}

\end{theorem}

\begin{proof}[Proof of Theorem 2]
Suppose $W$ has an inner product, which we use to identify $W$ with $V=W^*$.  Let $\Psi$ be the corresponding complement map, given by $\Psi(L)=\lin(L)$ as in Example \ref{ex.ipcm}. The relations coming out of Construction 1 then have the form $D_G(l_w-w)$ where $G$ is a face of $L$  and $w\in \lin(G)$. If $G$ is a one-dimensional cone, say $G=\Cone(w_i)$, then these relations have the form $D_i (l_{w_i}-w_i) $.  To see that these relations generate $D_G(l_w-w)$ for an arbitrary face $G$, simply note that $w\in \lin(G)$ implies that $w$ is a linear combination of the $w_i$ that generate $G$. Thus we see that $D_G(l_w-w)$ is in the ideal $J(L)$ defined in Step (e) of Construction 2.   In this way, we arrived exactly the statement of Construction 2. 
\end{proof}  

We have thus shown that Theorem 2 is a special case of Theorem 1. One could also prove Theorem 2 directly using an argument similar to the one we gave for Theorem 1.  Though we will not do this, we provide the following guidelines for how this would go.  These remarks are also useful for doing computations of $\mu$ in the inner product case.

The first step in proving Theorem 1 was the uniqueness of squarefree expressions.

\begin{proposition}\lbl{prop.unisf}  Let $W$ be an $n$-dimensional real vector space equipped with an inner product, and let $L$ be the basic cone in $W$ generated by the linearly independent set $\{w_1, \dots, w_k\}$ of $W$. Let $R(L)$ and $Z(L)$ be the rings defined in Construction \ref{constr.ip}. 
  Then squarefree monomials $\{D_G\ |\ G < L\}$ form a $\Lambda$-basis of $Z(L)$.
That is, every $\alpha\in Z(L)$ has a unique squarefree expression of the form
$$
\alpha = \sum_{G<L} \alpha_G D_G,
$$
with each $\alpha_G\in\Lambda$.
\end{proposition}

The proof is nearly identical to that of Proposition \ref{prop.unisfcm}. Instead of repeating the entire argument, we highlight the analogues of Lemma \ref{lem.standalone} and Equation (\ref{eq.didscm}), which can be used algorithmically to find squarefree expressions.

\begin{lemma}\lbl{lem.standaloneip}  Let $L=\Cone(w_1,\dots,w_k)$ be a basic cone.  Then
given $S\subset [k]$ and $i\in S$, there exists  $w\in W$ in the span of $\{w_j | j\in S\}$ such that  $[ w,  w_i] \neq 0$, but $[ w,w_j ]  = 0 $ for $j\in S\setminus \{i\}$. 
\end{lemma}

\begin{proof}
Define a map $f$ from the span of $\{w_j | j\in S\}$ to $\BR^{|S|}$ by $f(w) = ([ w,  w_j])_{j\in S}$.  It is easy to see that $f$ is linear and injective, hence surjective. The lemma follows.
\end{proof}

The existence of squarefree expressions boils down to computing a product $D_iD_S$, when $i\in S$. The key equation before was Equation (\ref{eq.didscm}). In the inner product case, take $w$ as in Lemma \ref{lem.standaloneip} and observe that $l_w - w$ is a linear combination of the  $l_{w_j} - w_j$ with $j\in S$.  Thus, since   $D_j (l_{w_j}-w_j) \in J$, it follows that $D_S (l_w- w) \in J$.  Thus in our ring $Z(L)$, we have the relation  
\begin{equation}\lbl{eq.dids}
 [ w,  w_i] D_i D_S + \sum_{j\notin S} [ w,  w_j] D_j D_S - w D_S = 0,
\end{equation}
which gives the desired squarefree expression for  $D_i D_S$ in $Z(L)$.

 The rest of the proof Proposition \ref{prop.unisf}  is nearly identical to the proof of Proposition \ref{prop.unisfcm}. As before the proof yields an algorithm for finding the squarefree expression for any polynomial in $R(L)$, which amounts to repeatedly using Equation (\ref{eq.dids}) with appropriately chosen $w$.

Compatability with inclusions in the inner product case is somewhat easier than in the complement map case.  

\begin{proposition}\label{prop.compatip}  Let $W$ be an $n$-dimensional vector space with an inner product, let $L$ be a basic cone in $N$, and let $\mu(L)\in \Lambda$ be as defined in Construction \ref{constr.ip}.  Suppose that $L$ lies in a rational subspace $W'\subset W$. Restrict the inner product on $W$ to $W'$, and let $\mu'(L)\in \Lambda'=\BC[[W']]$ result from Construction \ref{constr.ip}.  Then under the natural inclusion of $\Lambda'$ in $\Lambda$, the image of $\mu'(L)$ is $\mu(L)$.
\end{proposition} 
\begin{proof}
The inclusion of $\Lambda'$ in $\Lambda$ gives an inclusion $R'(L)=\Lambda'[[D_1, \dots, D_k]]\hookrightarrow R(L)=\Lambda[[D_1, \dots, D_k]]$.  It is easy to check that under this inclusion the ideal $J'(L)\subset J(L)$. By construction in $W'$, we have
$$
\Td(L) = \sum_{G<L} \lambda_G D_G +\gamma.
$$
with $\lambda_G \in \Lambda$ and $\gamma\in J'(L)$.  But then $\gamma\in J(L)$, so the above equation gives a squarefree expression for $\Td(L)$ in $Z(L)$. Uniqueness of squarefree expressions now forces $\mu(G)=\mu'(G)$ for all $G$.
\end{proof}

The remainder of the proof of Theorem 2 is then identical to the proof of Theorem 1 given in Section \ref{sec.ipproof}. One can simply ignore references to $\Psi$-genericity, and at the point when  Proposition \ref{prop.compatcm} was used, we instead deploy Proposition \ref{prop.compatip}.

\section{Examples and explicit formula}
In this section we present examples and calculations of values of $\mu$, leading up to an explicit formula for $\mu(L)$, where $L$ is a basic cone.\

\bigskip

{\bf 7.2 Calculation of $\mu$ in dimensions $0$, $1$, and $2$.}\ \ 
\begin{example}\lbl{ex.dim01}{\it (Dimensions 0 and 1.)} 
We fix a complement map $\Psi$ on the $n$-dimensional vector space $W$.  For the cone $L=\{0\}$, then it is easy to see that
$$
\mu(\{0\}) = 1,
$$
regardless of the complement map.  Now consider a one-dimensional basic cone $L=\Cone(w_1)$ that is generic with respect to $\Psi$. Let $u_1$ be the unique element of $\Psi(L)$ such that $\langle  w_1,u_1 \rangle=1$. The existence of $u_i$ follows from Lemma \ref{lem.standalone}. We will show that 
$$
\mu(L) = T(u_1) 
$$
where $T$ is the power series given by 
$$
T(z) = \frac{\td(z) - 1}{z} = \frac12 + \frac1{12} z -\frac1{720} z^3 + \cdots
$$
To see this, write
$$
\td(D_1) = 1 + D_1 T(D_1),
$$
and use the relation $D_1(l_{u_1}- u_1) = D_1( D_1 -u_1)$. Thus we can replace $D_1 T(D_1)$ with $D_1 T(u_1)$. We then find that $\mu(L)$, which is defined as the coefficient of $D_1$ in this expression, is equal to $T(u_1)$.
\end{example}

For cones of dimension two, we have the following.
\begin{proposition}\lbl{prop.explicit2d}
  Let $L=\Cone(w_1, w_2)$ be a two-dimensional basic dimensional cone in $W$ that is generic with respect to the complement map $\Psi$. For $i=1,2$, let  $u_i\in \Psi(\Span\{w_i\})$ satisfy  $\langle  w_i,u_i \rangle=1$.  For $i=1,2$, let $v_i$ in $\Psi(L)$  satisfy $\langle w_j,v_i \rangle=\delta_{i,j}$. Then 
  $$
  \mu(L) = T(v_1) T(v_2) + \frac{T(v_1)-T(u_1)}{v_2} + \frac{T(v_2)-T(u_2)}{v_1}.
  $$
  Furthermore the constant term is given by
  $$
  \mu_0(L)= \frac14 -\frac1{12} \bigl( \langle  w_2,u_1 \rangle + \langle  w_1,u_2 \rangle \bigr)   
  $$
\end{proposition}
\begin{proof}
Write 
$$
\td(D_1, D_2) = (1 +D_1T(D_1)) (1 + D_2T(D_2)) = 1 + D_1 T(D_1) + D_2 T(D_2) + D_1 D_2 T(D_1) T(D_2)
$$
We wish to compute the coefficient of $D_1 D_2$ in this expression.  Consider first the final term. Since $v_1\in\Psi(L)$, we have access to the relation
$D_1 D_2 (l_{v_1} - v_1 ) = D_1 D_2 ( D_1 - v_1 )= 0 $, and similarly $D_1 D_2 (D_2 - v_2)= 0$.   Hence the last term can be replaced by $D_1 D_2 T(v_1) T(v_2)$. 

Next we analyze the term $ D_1 T(D_1)$.  Here we use the relations $D_1 (l_{u_1} - u_1) =  D_1 ( D_1 + \langle  w_2,u_1 \rangle D_2 - u_1) $ and $D_1 D_2 (D_2 - v_2)$ to obtain
\begin{align*}
    D_1 T(D_1) &= D_1 T(-\langle  w_2,u_1 \rangle D_2 + u_1)\\
    &= D_1 T(u_1) + D_1 D_2  \biggr[ \frac{T(-\langle  w_2,u_1 \rangle D_2+ u_1)- T(u_1)}{D_2}  \biggl] \\
    &=D_1 T(u_1) + D_1 D_2 \biggr[ \frac{T(-\langle  w_2,u_1 \rangle v_2+ u_1)- T(u_1)}{v_2} \biggl]\\
    &= D_1 T(u_1) + D_1 D_2 \biggl[ \frac{T(v_1)-T(u_1)}{v_2} \biggr]
\end{align*}

Note that the bracketed expressions are actually power series in spite of the denominators.
In the final equality, we use the identity $ u_1 = \langle  w_1,u_1 \rangle v_1 + \langle  w_2,u_1 \rangle v_2$, which one can see by paring with $w_1$ and $w_2$.  

To summarize, the term $D_1 T(D_1)$ contributes to $\mu(L)$ a coefficient of $\frac{T(v_1)-T(u_1)}{v_2}$.  Similarly, $D_2 T(D_2)$ contributes $\frac{T(v_2)-T(u_2)}{v_1}$, and we are done.

As for the constant term $\mu_0(L)$, note that $T$ has constant term $1/2$, so $T(v_1) T(v_2) $ has constant term $1/4$.  For the other two terms, observe that
$$
T_2(z_1, z_2) := \frac{T(z_1+z_2)-T(z_1)}{z_2}
$$
is a power series in $z_1, z_2$ whose constant term equals the linear coefficient in the power series for $T(z)$, namely $1/12$.  It follows that $\frac{T(v_1)-T(u_1)}{v_2}$ has constant term $-\frac{1}{12} \langle  w_2,u_1 \rangle$.  A similar calculation applies to the remaining term, and the formula is established.

We note that two-dimensional $\mu_0$ values above can also be computed more directly in the spirit of [PT], where these $\mu_0$ were originally introduced.  To do so, we make the same computation but set all elements of $M$ to zero from the start.  That is, we work in the ring $Z_0(L)=Z(L)/A$  where all $A$ is the ideal generated by all $v\in M$. Now all of the relations in $J$ are homogeneous in the $D_i$, so to compute the coefficient of $D_1D_2$ in $\td(D_1, D_2)$, we must only consider the homogeneous terms of degree 2, which amounts to
$$\frac14 D_1 D_2 + \frac1{12}D_1^2+ \frac1{12} D_2^2= \biggl( \frac14 -  \frac1{12}\langle  w_2,u_1 \rangle - \frac1{12}\langle  w_1,u_2 \rangle \biggr) D_1 D_2,$$
 where we have used the relations   $D_1 l_{u_1}  =  D_1 ( D_1 + \langle  w_2,u_1 \rangle D_2 ) $ and $D_2 l_{u_2}  =  D_2 ( D_2 + \langle  w_1,u_2 \rangle D_1 ) $, which are valid in $Z_0(L)$.  Thus we see that $\mu_0(L)$, the coefficient of $D_1D_2$ in this expression, is as asserted in the proposition.
\end{proof}

\begin{example}\lbl{ex.2dip}
Specializing  Proposition \ref{prop.explicit2d}, we give the value of $\mu$ of a two-dimensional cone in the case where $\Psi$ comes from an inner product on $V$. 
Suppose $L=\Cone(w_1, w_2)$ is a two-dimensional basic dimensional cone in a two-dimensional space $W$ with a chosen inner product.   Let $\{v_1, v_2\}$ be the basis of $W$  dual to $\{w_1, w_2\}$  with respect to the inner product, so that   $[ v_i, w_j]=\delta_{i,j}$.  Note that the $u_i$ in the proposition are given by $u_i=\frac{w_i}{[w_i, w_i]}$.  Thus we obtain
  $$
  \mu(L) = T(v_1) T(v_2) + \frac1{v_2} \bigl({T(v_1)-T(\frac{w_1}{[w_1, w_1]})}\bigr) + \frac1{v_1} \bigl(T(v_2)-T(\frac{w_2}{[w_2, w_2]})\bigr), 
  $$
  with constant term
  \begin{equation}\lbl{eq.2dconst}
       \mu_0(L)= \frac14 -\frac1{12} \biggl( \frac{[w_1, w_2]}{[w_1, w_1]}  + \frac{[w_1, w_2]}{[w_2, w_2]}\biggr).    
  \end{equation}
\end{example}

\begin{example}\lbl{ex.2dtriangle}
We give a simple example of how these calculations of $\mu_0$ can be used to find the number of lattice points in a polytope. 
Consider the lattice $M=\BZ^2$ in  $V=\BR^2$ with the usual inner product, and let $P$ be the triangle with vertices $v_0=(0,0)$, $v_1=(t,0)$, and  $v_2=(0,t)$, where $t$ is a positive integer. The normal rays are $w_1=(0,1)$, $w_2=(1,0)$, and $w_0=(-1,-1)$.  Using Equation (\ref{eq.2dconst}), one computes $\mu_0$ of the $2$-dimensional normal cones as $\mu_0(\Cone(w_0, w_1))=\mu_0(\Cone(w_0, w_2))=\frac38$, and $\mu_0(\Cone(w_1, w_2))=\frac14$.
Using Example \ref{ex.dim01}, we see that $\mu_0(\Cone({w_i}))=\frac12$ for three one-dimensional cones, and $\mu_0(\{0\})=1$. Noting that the volume of $P$ is $t^2/2$,  the lattice-normalized volumes of all one dimensional faces equals $t$, and the volume of each zero-dimensional face is $1$, we obtain

\begin{align*}
    |P\cap \BZ^2| &= \sum_{F<P} \mu_0(C(P,F)) \Vol(F) \\ &= 1 \cdot \frac{t^2}{2} + 3\cdot \frac12 \cdot t + \frac38 \cdot 1 + \frac38 \cdot 1 + \frac14 \cdot 1 \\
    &=\frac12 t^2 +\frac32 t + 1.
\end{align*}

\end{example}

\bigskip\bigskip
{\bf 7.2 Formula for $\mu$ in arbitrary dimension.}\ \ 
We now extend these computations to basic cones of arbirtary dimension. 
Let $\Psi$ be a complement map and let $L=\Cone(w_1, \dots, w_k)$ be a $k$-dimensional cone that is generic with respect to $\Psi$.  Let $\Lambda$, $R(L)$, $J(L)$, $Z(L)$ be as in Construction \ref{constr.cm}.

For any subset $S\subset [k]$,  and $s\in S$, let $u_{S,s}$ denote the unique element in $\Psi(S)$ such that 
\begin{align}\lbl{eq.uequation}
    \langle w_s,u_{S,s} \rangle &= 1\ \ \ \ \text{and\ } \\
    \langle w_x,u_{S,s} \rangle &= 0\ \ \ \ \text{for all\ }x\in S\setminus \{s\}. 
\end{align}
Such $u_{S,s}$ exists and is unique by the argument of Lemma \ref{lem.standalone}.   Note that in the inner product case, in which we identify $V$ with $W$, we may take interpret the above equations as saying that $u_{S,s}$ is the unique element of $\Span \{w_s\ | \ s\in S \}$ such that 
\begin{align}\lbl{eq.u2equation}
    [ u_{S,s}, w_s ] &= 1\ \ \ \ \text{and\ } \\
    [ u_{S,s}, w_x] &= 0\ \ \ \ \text{for all\ }x\in S\setminus \{s\}. 
\end{align}

Further, if $s\notin S$,  we take $u_{S,s}=0$.
For any subset $T\subset S\subset[k]$, define $p_{S,T}$ (with $p$ standing for ``product"), by
$$
p_{S,T}= \prod_{t\in T} u_{S,t}, 
$$
 and let
$$
U_{S,T}= \frac{1}{p_{T,T}}\sum_{\mathcal{C}} (-1)^r \prod_{i=1}^{r} 
\frac1{p_{C_{i}, C_{i}\setminus C_{i-1}}}, 
$$
where the sum is taken over all chains $\mathcal{C}=(C_0, C_1, \dots, C_r)$  such that
$$
T=C_0\subsetneq C_1 \subsetneq \cdots \subsetneq C_r = S.
$$

For any $T\subset [k]$, let $u_T$ denote the vector $(u_{T,1}, \dots, u_{T,k}) $, remembering  that $u_{T,i}=0$ when $i\notin T$. For any polynomial $Q(D_1, \dots, D_k)$, we may form $Q(u_T)\in \Lambda$, by evaluation of this polynomial at the vector $u_T$. 
With this notation, we then have
\begin{theorem}\lbl{thm.explicit}   Fix a complement map $\Psi$ on $W$ and let $L$ be a $\Psi$-generic pointed cone in $W$. 
For any polynomial $Q(D_1, \dots, D_k)$, the squarefree expression of $Q$ in $Z(L)$ is given by
\begin{equation}\lbl{eq.explicit}
Q(D_1, \dots, D_k) = \sum_{S\subset [k]} \biggl\{
\sum_{T\subset S} Q(u_T) U_{S,T}
\biggr\} D_S
\end{equation}
\end{theorem}

\begin{proof} 
We fix the space $W$ of dimension $n$, and proceed by induction on $k$, the dimension of $L$. Note that if $k=0$, then we may take $Q$ to be a constant and both sides of our equation are seen to evaluate to this constant $Q$. 

Now suppose $k>0$.
It suffices to treat the case in which $Q$ is a monomial.
For $i=1, \dots, k$, let $v_i= u_{[k],i}$. We will apply the homomorphism  $e$ which sends $D_i$ to $v_i$.  Our first claim is that $e$ maps the left hand and right hand sides of our equation to the same element of $\Lambda$.   We will then show that this implies the equality.

When we apply $e$ to the left hand side, we get simply $Q(u_{[k]})$.  On the right hand side, noting that each $D_S$ maps to $p_{[k],S}$, we get a combination of the $Q(u_T)$, where the coefficient of each $Q(u_T)$ is given by
\begin{equation}\lbl{eq.uu}
    \sum_{S:T\subset S\subset [k]} p_{[k],S} U_{S,T}
\end{equation}
So it is enough to show that this expression equals 1 when $T=[k]$ and 0 otherwise. The first assertion is easy, since $p_{[k],[k]}U_{[k],[k]}=1$.  Suppose now that $T\subsetneq [k]$. Since $U_{S,T}$ is a sum over chains that start at $T$ and end at $S$, the entire epression in (\ref{eq.uu}) may be viewed as a sum over all chains that start at $T$ and end at an arbitrary subset $S$ such that $T\subset S\subset [k]$.  Pair each chain $\mathcal{C}=(C_0, C_1, \dots, C_r)$ for which the final set $S=C_r$ is not equal to $[k]$ with the chain $\hat{\mathcal{C}}=(C_0, C_1, \dots, C_r, [k])$ for which the final set is equal to $[k]$.  The contributions of $\mathcal{C}$  and $\hat{\mathcal{C}}$ to the sum in (\ref{eq.uu}) add up to 
\begin{equation}
    p_{[k],S}\frac{1}{p_{T,T}} (-1)^r \prod_{i=1}^{r} \frac1{p_{C_{i}, C_{i}\setminus C_{i-1}}}
    +
    p_{[k],[k]}\frac{1}{p_{T,T}} (-1)^{r+1} \prod_{i=1}^{r} \frac1{p_{C_{i}, C_{i}\setminus C_{i-1}}}\cdot 
    \frac1{p_{[k],[k]\setminus S}}
    =0.
\end{equation}
Hence, the entire sum in (\ref{eq.uu}) vanishes, as desired.

Thus we have established that the two sides of the equation (\ref{eq.explicit}) have the same image under $e$.  We now establish the equality itself.  For any $S\subset [k]$, let $\beta_S$ be the bracket-enclosed coefficient of $D_S$ on the right hand side of Equation (\ref{eq.explicit}).  We wish to show that $\beta_S=\alpha_S$, where $\alpha_S$ is the correct coefficient in the squarefree expression
\begin{equation}\lbl{eq.qda}
Q(D_1, \dots, D_k) = \sum_{S\subset [k]} \alpha_S D_S .
\end{equation}
We first establish that $\beta_S=\alpha_S$ for any $S\subsetneq [k]$. We can assume $S=[j]$, with $j<k$. If $\Supp(Q)\not\subset S$, then we have $\alpha_S=0$ by Lemma \ref{lem.suppgrow}.  But also $\beta_S=0$.  This follows since for $T\subset S$, $\Supp(Q)\not\subset T$. So taking $i\in \Supp(Q)\setminus T$, we have $u_{T,i}=0$, whence $Q(u_T)=0$.

Thus we may suppose that $\Supp(Q)\subset S$, so $Q$ is a monomial in $D_1, \dots D_j$, and since $j<k$, and our induction hypothesis asserts that 
\begin{equation}\lbl{eq.Qalpha}
Q(D_1, \dots, D_k) =\sum_{S\subset [j]} \beta_S D_S
\end{equation}
in the ring $Z(L')$, where $L'$ is the face of $L$, spanned by $w_1, \dots, w_j$. Lemma \ref{lem.facesf} allows us to conclude that $\alpha_S=\beta_S$.  

Having established that $\alpha_S=\beta_S$ for all $S\subsetneq [k]$, we must only establish that  $\alpha_{[k]}=\beta_{[k]}$.  But 
$$
\sum e(D_S) \alpha_S = e(Q(D_1, \dots, D_j)) = \sum e(D_S) \beta_S,
$$
the first equation following from Equation (\ref{eq.qda}), and the second from our previously-established assertion that $e$ preserves our identity. Since $e(D_{[k]})\neq 0$ and $\alpha_S=\beta_S$ for all $S\subsetneq [k]$, this forces $\alpha_{[k]}=\beta_{[k]}$.
\end{proof}

Taking $Q$ to be the Todd polynomial, we obtain
\begin{theorem}\lbl{thm.explicittd}  With the above set-up,
$$
\mu(L) = \sum_{T\subset [k]} \td(u_T) U_{[k],T}
$$

\end{theorem}

We remark that Diaz, Le, and Robins, in their study of solid angle sums \cite{DLR}, uncovered ingredients that resemble our chain sums $U_{S,T}$, at least superficially.  It would be interesting to look for deeper connections here. 

\bigskip
\bigskip We illustrate Theorems \ref{thm.explicit} and \ref{thm.explicittd} in the case of two-dimensional cones.   Let $L=\Cone(w_1, w_2)$ be a two-dimensional basic cone in the $n$-dimensional space $W$.  We must compute $U_{[2], T}$, for each subset of $T\subset [2]$.  For $T=[2]$, there is only one chain, with $r=0$, yielding
$$
U_{12,12}= \frac1{u_{12,1}u_{12,2}}.
$$
Here we write $U_{12,12}$ in place of the more cumbersome $U_{\{1,2\},\{1,2\}}$, etc. 
For $T=\{1\}$, there is again only one chain, this time with $r=1$, and we see that
$$
U_{12,1}= -\frac1{u_{1,1}u_{12,2}}.
$$
The situation is similar for $T=\{2\}$.
For $T=\emptyset$, there are now three chains, $\emptyset \subset \{1, 2\}$  with $r=1$, and two chains of the form $\emptyset \subset \{i\}\subset \{1, 2\}$ with $r=2$.  One obtains
$$
U_{12,\emptyset}= -\frac1{u_{12,1}u_{12,2}}+\frac1{u_{12,2}u_{1,1}} +\frac1{u_{12,1}u_{2,2}}.
$$
One thus obtains
\begin{equation}\lbl{eq.2d}
\mu(L) = \frac{\td(u_{12,1},u_{12,2})}{u_{12,1}u_{12,2}}- \frac{\td(u_{1,1})}{u_{1,1}u_{12,2}}
- \frac{\td(u_{2,2})}{u_{2,2}u_{12,1}}
-\frac1{u_{12,1}u_{12,2}}+\frac1{u_{12,2}u_{1,1}} +\frac1{u_{12,1}u_{2,2}}
\end{equation}
In the notation of Proposition \ref{prop.explicit2d}, $u_{12,i}$ is denoted there by $v_i$, and $u_{i,i}$ is denoted there $u_i$.  One then checks agreement of Equation (\ref{eq.2d}) with Proposition \ref{prop.explicit2d} using a straightforward calculation.

\section{Projective Space Example}

In this section, we consider the case where our polytope $P$ is the standard $n$-simplex in $\BR^{n+1}$, whose corresponding toric variety is the projective space $\BP^n$. The work of Diaconis and Fulton [DF] gives a nice formula for the (ordinary) Todd class of projective space.  Building upon this work, we begin a corresponding analysis of the equivariant Todd class of projective space, and the closely related problem of finding an Euler-Maclaurin formula for the simplex.  In particular, the approach of [DF] suggests a certain complement map, and we analyze the corresponding $\mu$ values arising out of our Construction \ref{constr.cm}.  Interestingly, this complement map arises neither from an inner product nor from a flag. 

First let $N$ be a full-rank lattice in the $n$-dimensional vector space $W$. Let $\{w_1, \dots, w_n\} $ be a basis of the lattice $N$ and let $w_0=-\sum_{i=1}^n w_i$.  Then the rays $w_0, \dots, w_n$ determine a fan in $N$ with $n+1$ maximal cones, all of which are basic.  This is the fan whose toric variety is the projective space $\BP^n$. 
This fan can be realized as the inner normal fan of a simplex.  Specifically, if $M$ is the dual lattice in the dual space $V$, with $\{v_1, \dots, v_n\}$ the basis dual to $\{w_1, \dots, w_n\}$, then we may take the simplex $P$ to be the convex hull of $\{0,v_1, \dots, v_n\} $. 


 Diaconis and Fulton found an expression for the (ordinary) Todd class of $\BP^n$.  Their computation amounts to choosing the following complement map on $W$:  For any $i=0, \dots, n$, let $\Psi(\Cone(w_i))= \Span \{u_i\} $, where $u_i$ is the element of $V$ that evaluates to $1$ on $w_i$ and $-1$ on $w_{i+1}$, where indices are taken cyclically modulo $n+1$. That is, we set $u_i=v_i - v_{i+1}$ for $i=1,...,n-1$, $u_n=v_n$ and  $u_0 = - v_1$. 
We wish to find expressions for the $\mu$ values of the cones of this fan.  By symmetry, we  choose the cone $L=\Cone(w_1, \dots,w_n)$ and consider the corresponding ring
 $$
 Z(L)=\Lambda[[D_1, \dots, D_n]]/J(L),
 $$
where 
$$
J(L)=\biggl< D_i(D_i- D_{i+1}- u_i) \ | \ i=1,\dots n \biggr>
$$
where we interpret $D_{n+1}$ as $0$. 
(We remark that equivalent computaions could be carried out in the global ring $Z(\Sigma)$ mentioned in Section 1.6. However, in the spirit of our constructions, we have chosen to work locally in $Z(L)$.)

Taking all $u_i=0$ specializes the ring $Z(L)$ to the ring $Z_0(L)= \BC[[D_1\dots D_n]]/J_0(L)$, where $J_0(L)$ is the ideal,
$$
J_0(L)=\biggl< D_i(D_i- D_{i+1}) \ | \ i=1,\dots n \biggr>. 
$$
This is the ring considered in [DF, Sections 5 and 6], which is used to give a geometric expression for the ordinary Todd class of projective space.  Place the fan defining projective space $\BP^n$ symmetrically in $\BR^{n+1}$ using  the $n+1$  rays $\rho_i=e_i-e_{i+1}$, $i=1, \dots, n+1$, where the $e_i$ are the standard basis vectors in $\BR^{n+1}$ indexed cyclically.  The Diaconis-Fulton evaluation of the Todd coefficients amounts to the formula in $J_0(L)$
\begin{equation}
    \td(D_1,\dots, D_n) = \sum_{S\subset \{1, \dots, n\}} f_S D_S,
\end{equation}
 where $f_S$ denotes the fraction of the linear space spanned by the vectors $\{\rho_s | s\in S\} $ that is occupied by the cone over these vectors. Given such a nice geometric interpretation of the ordinary Todd class coefficients, it would be wonderful if in the equivariant case, there were also a geometric way to express the Todd coefficients coming out of this same complement map. Though we do not know how to do this in general, we offer the following calculations for the cones of $\Sigma$ of dimensions $1$ and $2$.

 


Consider the cone $L=\Cone(w_1, \dots, w_n)$ with the complement map $\Psi$ described above.  For the one-dimensional face $L'=\Cone(w_i)$
we have the ideal $J(L')=D_i(D_i-u_i)$, and we see as in Example \ref{ex.dim01} that
\begin{equation}
    \mu(L') = T(u_i),
\end{equation}
  The constant term is $1/2$, which is indeed the fraction of the line through $\rho_i=e_i-e_{i+1}$ that is occupied by ray through $\rho_i$. 

Now consider a two-dimensional face of $L$, say $L'=\Cone(w_i, w_j)$.  Now there are distinct cases depending on whether $i$ and $j$ are consecutive integers or not. First suppose that $i$ and $j$ are not consecutive.  We work modulo the ideal $J(L')$ generated by $D_i(D_i-u_i)$ and $D_j(D_j-u_j)$.   Calculating directly or invoking Proposition \ref{prop.explicit2d}, we see that 
\begin{equation}
    \mu(L') = T(u_i) T(u_j).
\end{equation}
The constant term here is $1/4$, which indeed is the fraction the linear span of $\rho_i=e_i-e_{i+1}$ and $\rho_j=e_j-e_{j+1}$ occupied by the cone generated by these two vectors. 

Now suppose $i$ and $j$ are consecutive, say $j={i+1}$. This time the ideal 
 $J(L')$ is generated by $D_i(D_i-D_{i+1}-u_i)$ and $D_{i+1}(D_{i+1}-u_{i+1})$.  Now we find
\begin{equation}
    \mu(L') = T_2(u_i, u_{i+1}) + T(u_i+u_{i+1})T(u_{i+1}).
\end{equation}
Here the constant term is $1/12 + 1/4 =1/3$, which equals the fraction of the linear span of $\rho_i$ and $\rho_{i+1}$ taken up by the cone generated by these two vectors.

\bigskip

\address{{\tt bfischer@seattleacademy.org}:  Seattle Academy, Seattle WA.}

\address{{\tt jamie@reed.edu}: Reed College, Portland OR.}

\ifx\undefined\bysame
        \newcommand{\bysame}{\leavevmode\hbox
to3em{\hrulefill}\,}
\fi


\begin{thebibliography}{[EMSS]}

\bibitem[Bv]{Bv} A. Barvinok,
        {\em  A polynomial time algorithm for counting integral points 
        in polyhedra when the dimension is fixed},  
        Math. Oper. Res. {\bf 19}  (1994) 769--779.
 
\bibitem[BP]{BP} A. Barvinok and J.E. Pommersheim, 
        {\em An algorithmic theory of lattice points in polyhedra},
        Math. Sci. Res. Inst. Publ., {\bf 38} (1999) 91--147.

\bibitem[BV1]{BV1} N. Berline and M. Vergne,
        {\em Local Euler-Maclaurin formula for polytopes},  
        Mosc. Math. J.  {\bf 7}  (2007) 355--386.
        
\bibitem[BV2]{BV2} N. Berline and M. Vergne,
        {\em The equivariant Todd genus of a complete toric variety, with Danilov condition},  
        J. Alegbra  {\bf 313}  (2007) 28--37.

\bibitem[BR]{BR} M. Beck and S. Robins,
	{\em Computing the continuous discretely: integer-point enumeration in polyhedra},
	Springer-Verlag 2007.

\bibitem[Br1]{Br} M. Brion,
        {\em  Points entiers dans les poly\`edres convexes},  
        Ann. Sci. \'Ecole Norm. Sup. {\bf 4} 21 (1988) 653--663. 

\bibitem[Br2]{Br2} M. Brion,
        {\em  Equivariant Chow groups for torus actions},  
        Transformation groups {\bf 2} (1997) 225--267. 

\bibitem[BV3]{BV3} M. Brion and M. Vergne,
	    {\em An equivariant Riemann-Roch theorem for complete, simplicial toric varieties},
	    J. reine angew. Math. {\bf 482} 67--92.
	
\bibitem [CK]{CK} A. Connes and D. Kreimer,
        {\em Hopf algebras, Renormalization and Noncommutative Geometry},
        Comm. Math. Phys. {\bf 199} (1988) 203-242.

	
\bibitem[DF]{DF} P. Diaconis and W. Fulton, 
        {\em A growth model, a game, an algebra, Lagrange inversion, and characteristic classes}, Rend. Sem. Mat. Univ. Pol. Torino {\bf 49} (1991), 95–-119.
          
\bibitem[DLR]{DLR} R. Diaz, Q.-N Le, S. Robins, 
        {\em Fourier transforms of polytopes, solid angle sums, and discrete volume}, 2016 http://arxiv.org/abs/1602.08593.

        
\bibitem[Fi]{Fi} B. Fischer, 
        {\em Perturbed Polyhedra and the Construction of Local Euler-Maclaurin Formulas},
        Ph.D. thesis, Boston University, 2016.

\bibitem[FP]{FP} B. Fischer and J. Pommersheim,
        {\em Cycle-Level Products in Equivariant Cohomology of Toric Varieties}
        Michigan Math. J. {\bf 63: 4} (2014), 845--864.


\bibitem[Fu]{Fu} W. Fulton,
        {\em Introduction to toric varieties}.
        Annals of Mathematics Studies {\bf 131} Princeton University Press
        1993.
        
 \bibitem[Fu2]{Fu2} W. Fulton,
        {\em Equivariant Cohomology in Algebraic Geometry},
        Lecture Thirteen: Toric Varieties.  
        
\bibitem[GP]{GP} S. Garoufalidis and J. Pommersheim,
		{\em Sum-Integral interpolators and the Euler-Maclaurin formula for polytopes}.
		Trans. Amer. Math. Soc. {\bf 364} (2012) 2933--2958.
		
\bibitem[GPZ]{GPZ} L. Guo, S. Paycha and B. Zhang,
        {\em Algebraic Birkhoff Factorisation and the Euler-Maclaurin formula on cones}.
        arXiv:1306.3420.

\bibitem[KP]{KP} A. G. Khovanskii and A. V. Pukhlikov, 
	{\em A Riemann-Roch theorem for integrals and sums of quasipolynomials over virtual polytopes},
	St Petersburg Math. {\bf 4} (1993) 789--812.

\bibitem[La]{La} J. Lawrence, 
        {\em Rational-function-valued valuations on polyhedra},
        DIMACS Ser. Discrete Math. Theoret. Comput. Sci.,{\bf 6} 
        (1991) 199--208.

\bibitem[McM]{McM} P. McMullen, 
	  {\em Weakly continuous valuations on convex polytopes},
	  Archiv Math. {\bf 41} (1983) 555--564.

\bibitem[Mo]{Mo} R. Morelli,
        {\em Pick's theorem and the Todd class of a toric variety},  
        Adv. Math.  {\bf 100}  (1993) 183--231.

\bibitem[Pi]{Pi} G.A. Pick,
        {\em Geometrisches zur Zahlenlehre}, 
        Sitzenber. Lotos (Prague) {\bf 19} (1899) 311--319.

\bibitem[PT]{PT} J. Pommersheim and H. Thomas, 
        {\em Cycles representing the Todd class of a toric variety},
        J. Amer. Math. Soc. {\bf 17} (2004) 983--994.

\bibitem[Th]{Th} H. Thomas,
        {\em Cycle-level intersection theory for toric varieties}, 
        Canad. J. Math.  {\bf 56}  (2004) 1094--1120. 

\end{thebibliography}
\end{document}